\documentclass[]{preprint}
\usepackage{times}
\usepackage{mhequ}
\usepackage{mhfig}
\usepackage{mhenvs}
\usepackage{Symbols}
\usepackage{ScriptFonts}

\def\prob{{\bf P}}
\def\expect{{\bf E}}
\def\from{:}

{
\theorembodyfont{\rmfamily}
\newtheorem{notation}[lemma]{Notation}
}

\ifx\pdfoutput\undefined
\else
	\input protcode.tex
\fi

\makeatletter
\def\Kern{\@ifnextchar({\Kern@par}{\Kern@nopar}}
\def\Kern@par(#1){(#1)^{H-{1\over 2}}}
\def\Kern@nopar#1{{#1}^{H-{1\over 2}}}
\def\maketag#1#2{(#1)\def\@currentlabel{#1}\label{#2}}
\makeatother

\def\Noise{\CW}
\def\State{\CX}
\def\pNoise{{\prob_{\! w}}}
\def\Evol{{\cscr Q}}
\def\Coupl{{\cscr C}}
\def\BigSpace{{\cscr X}}
\def\Dens{{\cscr D}}
\def\C#1{{C^{\mathbf{A1}}_{#1}}}
\def\nb{\penalty10000}

\def\pMeas{{\cscr M}_1}
\def\Meas{{\cscr M}_\phi}
\def\cadlag{c\`adl\`ag}
\def\Wien{{\mathsf W}}
\def\TV{{\mathrm TV}}
\def\step#1{{\textbf #1}}

\begin{document}
\ifx\pdfoutput\undefined
\else
    \setprotcode\font
    {\it \setprotcode\font}
    {\bf \setprotcode\font}
    {\bf \it \setprotcode\font}
    \pdfprotrudechars=1
\fi

\title{Ergodicity of Stochastic Differential Equations Driven by
Fractional Brownian Motion}

\author{Martin Hairer}
\institute{Mathematics Research Centre, University of Warwick\\
\email{hairer@maths.warwick.ac.uk}\\
Homepage: \texttt{http://www.maths.warwick.ac.uk/\char`~hairer/}}

\maketitle

\begin{abstract}
We study the ergodic properties of finite-dimensional systems of SDEs driven by 
non-degenerate additive fractional Brownian motion with
arbitrary Hurst parameter $H\in(0,1)$. A general framework is constructed to make
precise the notions of ``invariant measure'' and ``stationary state'' for such a system.
We then prove under rather weak dissipativity conditions that such 
an SDE possesses a unique stationary solution and that the convergence
rate of an arbitrary solution towards the stationary one is (at least) algebraic. A lower bound on the
exponent is also given.
\end{abstract}

\tableofcontents

\section{Introduction and main result}

In  this paper, we investigate the long-time behaviour of stochastic differential equations driven by fractional
Brownian motion. Fractional Brownian motion (or fBm for short) is a centred Gaussian process satisfying $B_H(0) = 0$ and
\begin{equ}[e:deffBm]
 \expect|B_H(t) - B_H(s)|^2 = |t-s|^{2H}\;,\qquad t,s>0\;,
\end{equ}
where $H$, the Hurst parameter, is a real number in the range $H\in(0,1)$. When $H={1\over 2}$, one recovers 
of course the usual Brownian motion, so this is a natural one-parameter family of generalisations
of the ``standard'' Brownian motion. It follows from \eref{e:deffBm} that fBm is also self-similar, but with the scaling law
\begin{equ}
t \mapsto B_H(at) \qquad\approx\qquad t \mapsto a^H B_H(t)\;,
\end{equ}
where $\approx$ denotes equivalence in law. Also, the sample paths of $B_H$ are $\alpha$-H\"older continuous for every
$\alpha < H$. The main difference between fBm and the usual Brownian motion is that it is neither
Markovian, nor a semi-martingale, so most standard tools from stochastic calculus 
cannot be applied to its analysis.

Our main motivation is to tackle the problem of ergodicity in non-Markovian systems.
Such systems arise naturally in several situations. In physics, stochastic forces are used to describe the
interaction between a (small) system and its (large) environment. There is no a-priori reason to assume that
the forces applied by the environment to the system are independent over disjoint time intervals.
In statistical mechanics, for example, a  non-Markovian
noise term appears when one attempts to derive 
the Langevin equation from first principles \cite{MR98h:82049,MR90a:82002}. 
Self-similar stochastic processes like fractional Brownian motion
appear naturally in hydrodynamics \cite{MR39:3572}.
It appears that fractional Brownian motion is also useful to model long-time correlations in stock markets \cite{MR2001g:60129,fbmfinance}.

Little seems to be known about the long-time behaviour of non-Markovian systems.
In the case of the non-Markovian Langevin equation (which is \textit{not} covered 
by the results in this paper
due to the presence of a delay term), the stationary solution is explicitly known to be distributed according to
the usual equilibrium Gibbs measure. The relaxation towards equilibrium is a very hard problem that was
solved in \cite{MR98h:82049,MR2001c:82006}. It is however still open in the non-equilibrium case, where the
invariant state can not be guessed a-priori. 
One well-studied general framework for the study of systems driven by noise with extrinsic 
memory like the ones considered in this paper is given by the theory of Random Dynamical Systems (see the
monograph \cite{RDS} and the reference list therein). In that framework, the existence of random attractors,
and therefore the \textit{existence} of invariant measures
seems to be well-understood. On the other hand, the problem of \textit{uniqueness} (in an appropriate sense, see
the comment following \theo{theo:largeH} below) of the invariant measure on the
random attractor seems to be much harder, unless one can show that the system possesses a unique
stochastic fixed point. The latter situation was studied in
\cite{SchmalMasl} for infinite-dimensional evolution equations driven by fractional Brownian motion. 

The reasons for choosing fBm as driving process for \eref{e:mainequ} below are twofold. First, in particular
when $H > {1\over 2}$, fractional Brownian motion presents genuine long-time correlations that
persist even under rescaling. The second reason is that there exist simple, explicit formulae that relate
fractional Brownian motion to ``standard'' Brownian motion, which simplifies our analysis.
We will limit ourselves to the case where the memory of the system
comes entirely from the driving noise process, so we do not consider stochastic delay equations.

We will only consider equations driven by non-degenerate additive noise, \ie\ we consider equations of the form
\begin{equ}[e:mainequ]\tag{SDE}
dx_t = f(x_t)\,dt + \sigma\,dB_H(t)\;,\qquad x_0 \in \R^n\;,
\end{equ}
where $x_t \in \R^n$, $f\from \R^n \to \R^n$, $B_H$ is an $n$-dimensional fractional Brownian motion with Hurst parameter $H$,
and $\sigma$ is a constant and invertible $n\times n$ matrix. Of course, \eref{e:mainequ} should be interpreted as an integral
equation.

In order to ensure the existence of globally  bounded solutions and in order to have some control on the speed at which
trajectories separate, we make throughout the paper the following assumptions on the components of \eref{e:mainequ}:
\begin{claim}[A2']
\item[\textbf{A1}] \textit{Stability.} There exist constants $\C{i} > 0$ such that
\begin{equ}
\scal{f(x)-f(y),x-y} \le \min\{\C1 - \C2\|x-y\|^2, \C3\|x-y\|^2\}\;,
\end{equ}
for every $x,y \in \R^n$.
\item[\textbf{A2}] \textit{Growth and regularity.} There exist constants $C,N>0$ such that $f$ and its derivative satisfy
\begin{equ}
\|f(x)\| \le C\sml(1+\|x\|\smr)^N\;, \qquad \|Df(x)\| \le C\sml(1+\|x\|\smr)^N\;,
\end{equ}
for every $x \in \R^n$.
\item[\textbf{A3}] \textit{Non-degeneracy.} The $n\times n$ matrix $\sigma$ is invertible. 
\end{claim}

\begin{remark}\label{rem:sigma}
We can assume that $\|\sigma\| \le 1$ without any loss of generality. This assumption 
will be made throughout the paper in order to simplify some expressions.
\end{remark}

One typical example that we have in mind is given by
\begin{equ}
f(x) = x - x^3\:,\qquad x\in \R\;,
\end{equ}
or any polynomial of odd degree with negative leading coefficient.
Notice that $f$ satisfies \textbf{A1}--\textbf{A2}, but that it is not globally Lipschitz continuous.

When the Hurst parameter $H$ of the fBm driving \eref{e:mainequ} is bigger than $1/2$, 
more regularity for $f$ is required, and we will then sometimes assume that the following 
stronger condition holds instead of \textbf{A2}:

\begin{claim}[A2']
\item[\textbf{A2'}] \textit{Strong regularity.} The derivative of $f$ is globally bounded.
\end{claim}

Our main result is that \eref{e:mainequ} possesses a \textit{unique} stationary solution.
Furthermore, we obtain an explicit bound showing that every (adapted) solution to 
\eref{e:mainequ} converges towards this stationary solution, and that this convergence
is at least algebraic. We make no claim concerning the optimality of this bound for the class
of systems under consideration.
Our results are slightly different for small and for large values of $H$, so we state them separately.

\begin{theorem}[Small Hurst parameter]\label{theo:smallH}
Let $H\in \sml(0,{\textstyle{1\over 2}}\smr)$ and let $f$ and $\sigma$ satisfy \textbf{A1}--\textbf{A3}.
Then, for every initial condition, the solution to \eref{e:mainequ} converges towards a unique stationary solution
in the total variation norm. Furthermore, for every $\gamma < \max_{\alpha < H} \alpha(1-2\alpha)$, the difference between the solution and the stationary solution
is bounded by $C_\gamma t^{-\gamma}$ for large $t$.
\end{theorem}

\begin{theorem}[Large Hurst parameter]\label{theo:largeH}
Let $H\in \sml({\textstyle{1\over 2}},1\smr)$ and let $f$ and $\sigma$ satisfy \textbf{A1}--\textbf{A3} and \textbf{A2'}.
Then, for every initial condition, the solution to \eref{e:mainequ} converges towards a unique stationary solution
in the total variation norm. Furthermore, for every $\gamma < {1\over 8}$, the difference between the solution and the stationary solution
is bounded by $C_\gamma t^{-\gamma}$ for large $t$.
\end{theorem}

\begin{remark}
The ``uniqueness'' part of these statements should be understood as uni\-que\-ness in law in 
the class of stationary solutions
adapted to the natural filtration induced by the two-sided fBm that drives the equation. There could in theory be other stationary solutions,
but they would require knowledge of the future to determine the present, so they are usually discarded as unphysical.

Even in the context of Markov processes, similar situations do occur. One can well have uniqueness of the invariant measure,
but non-uniqueness of the stationary state, although other stationary states would have to foresee the future. In this sense,
the notion of uniqueness appearing in the above statements is similar to the notion of uniqueness of the invariant measure for Markov processes. (See 
\eg \cite{RDS}, \cite{MR93c:60096} and \cite{MR1888416}
for discussions on invariant measures that are not necessarily measurable with respect to the past.)
\end{remark}

\begin{remark}
The case $H={\textstyle{1\over 2}}$ is not covered by these two theorems, but it is well-known that the convergence 
toward the stationary state is exponential in this case (see for example \cite{MT}).
In both cases, the word  ``total variation'' refers to the total variation distance between measures on the space
of paths, see also \theo{theo:mainTheoFormal} below for a rigorous formulation of the results above.
\end{remark}

\subsection{Idea of proof and structure of the paper}

Our first task is to make precise the notions of ``initial condition'', ``invariant measure'', ``uniqueness'', and ``convergence''
appearing in the formulation of Theorems~\ref{theo:smallH} and \ref{theo:largeH}. This will be
achieved in Section~\ref{sec:general} below, where we construct a general framework for the study
of systems driven by non-Markovian noise. Section~\ref{sec:setting} shows how \eref{e:mainequ}
fits into that framework.

The main tool used in the proof of Theorems~\ref{theo:smallH} and \ref{theo:largeH} is a coupling
construction similar in spirit to the ones presented in \cite{MatNS,HExp02}. More precisely, we first
show by some compactness argument that there exists at least one invariant measure $\mu_*$ for
\eref{e:mainequ}. Then, given an initial condition distributed according to some arbitrary measure $\mu$, 
we construct a ``coupling process''
$(x_t,y_t)$ on $\R^n\times \R^n$ with the following properties:
\begin{claim}[3.]
\item[1.] The process $x_t$ is a solution to \eref{e:mainequ} with initial condition $\mu_*$.
\item[2.] The process $y_t$ is a solution to \eref{e:mainequ} with initial condition $\mu$.
\item[3.] The random time $\tau_\infty = \min\{t\,|\, x_s = y_s\quad \forall s \ge t\}$ is almost surely finite.
\end{claim}
The challenge is to introduce correlations between $x_s$ and $y_s$ in precisely such a way that
$\tau_\infty$ is finite. If this is possible, the uniqueness of the invariant measure follows immediately.
Bounds on the moments of $\tau_\infty$ furthermore translate into bounds on the rate of convergence towards this
invariant measure.
In Section~\ref{sec:coupl}, we expose the general mechanism by which
we construct this coupling. Section~\ref{sec:defCoupl} is then devoted to the precise formulation of the
coupling process and to the study of its properties, which will be used in Section~\ref{sec:mainproof}
to prove Theorems~\ref{theo:smallH} and \ref{theo:largeH}. We conclude this paper with a few
remarks on possible extensions of our results to situations that are not covered here.

\begin{acknowledge}
The author wishes to thank Dirk Bl\"omker, David Elworthy, Xue-Mei Li, Neil O'Connell, and 
Roger Tribe for their interest in this work and for many helpful suggestions, stimulating questions, and 
pertinent remarks. He would also like to thank the referee for his careful reading of the manuscript.

The author also wishes to thank the Mathematics Research Centre of the University of Warwick 
for its warm hospitality. This work was supported by the Fonds National Suisse.
\end{acknowledge}

\section{General theory of stochastic dynamical systems}
\label{sec:general}

In this section, we first construct an abstract framework that can be used to model a large
class of physically relevant models where the driving noise is stationary. Our framework
is very closely related to the framework of random dynamical systems with however one 
fundamental difference. In the theory of random dynamical systems (RDS), the abstract space
$\Omega$ used to model the noise part typically encodes the {\em future}
of the noise process. In our framework of ``stochastic dynamical systems'' (SDS)
the noise space $\Noise$ typically encodes the {\em past} of the noise process. As a consequence,
the evolution on $\Noise$ will be stochastic, as opposed to the deterministic evolution on $\Omega$
one encounters in the theory of RDS. This distinction may seem futile at first sight, and one could
argue that the difference between RDS and SDS is non-existent by adding the past of the noise process to
$\Omega$ and its future to $\CW$.

The additional structure we require is that the evolution on $\Noise$
possesses a \textit{unique} invariant measure. Although this requirement may sound very strong, it is actually
not, and most natural examples satisfy it, as long as $\Noise$ is chosen in such a way that it
does not contain information about the future of the noise. 
In very loose terms, this requirement of having a unique invariant 
measure states that the noise process driving our system is stationary and that the Markov process
modelling its evolution captures all its essential features in such a way that it could not be
used to describe a noise process different from the one at hand. In particular, this means that
there is a continuous inflow of ``new randomness'' into the system, which is a crucial feature
when trying to apply probabilistic methods to the study of ergodic properties of the system. 
This is in opposition to the RDS formalism, where the noise is ``frozen'', as soon as an 
element of $\Omega$ is chosen. 

From the mathematical point of view, we will consider that the physical process we are interested in
lives on a ``state space'' $\State$ and that its driving noise belongs to a 
``noise space'' $\Noise$. In both cases, we only consider Polish (\ie complete, separable, and
metrisable) spaces. One should think of the state space as a relatively small space which contains
all the information accessible to a physical observer of the process. The noise space should be 
thought of as a much bigger abstract space containing all the information needed to construct
a mathematical model of the driving noise up to a certain time. The information contained in the 
noise space is not accessible to the physical observer. 

Before we state our definition of a SDS, we will recall several notations and definitions, mainly
for the sake of mathematical rigour. The reader can safely skip the next subsection and come back
to it for reference concerning the notations and the mathematically precise definitions of the concepts
that are used.

\subsection{Preliminary definitions and notations}

First of all, recall he definition of a transition semigroup:
\begin{definition}\label{def:kernel}
Let $(\CE,{\cscr E})$ be a Polish space endowed with its
Borel $\sigma$-field. A {\em transition semigroup} $\CP_t$ on $\CE$ is a family of maps
 $\CP_t \from \CE \times {\cscr E} \to [0,1]$ indexed by $t \in [0,\infty)$ such that
\begin{claim}[iii)]
\item[i)] for every $x\in\CE$, the map $A \mapsto \CP_t(x,A)$ is a probability measure on $\CE$
and, for every $A \in {\cscr E}$, the map $x \mapsto \CP_t(x,A)$ is ${\cscr E}$-measurable,
\item[ii)] one has the identity
\begin{equ}
\CP_{s+t}(x,A) = \int_{\CE} \CP_s(y,A)\,\CP_t(x,dy)\;,
\end{equ}
for every $s,t > 0$, every $x\in\CE$, and every $A\in {\cscr E}$.
\item[iii)] $\CP_0(x,\cdot) = \delta_x$ for every $x\in\CE$. 
\end{claim}
\end{definition}
We will freely use the notations 
\begin{equ}
\sml(\CP_t \psi\smr)(x) = \int_\CE \psi(y)\,\CP_t(x,dy)\;,\qquad \sml(\CP_t \mu\smr)(A) = \int_\CE \CP_t(x,A)\,\mu(dx)\;,
\end{equ}
where $\psi$ is a measurable function on $\CE$ and $\mu$ is a measure on $\CE$.

Since we will always work with topological spaces, we will require our transition semigroups
to have good topological properties. Recall that a sequence $\{\mu_n\}$ of measures on a topological
space $\CE$ is said to converge toward a limiting measure $\mu$ in the weak topology if
\begin{equ}
\int_\CE \psi(x)\,\mu_n(dx) \to \int_\CE \psi(x)\,\mu(dx)\;,\qquad \forall \psi \in \CC_b(\CE)\;,
\end{equ}
where $\CC_b(\CE)$ denotes the space of bounded continuous functions from $\CE$ into $\R$.
In the sequel, we will use the notation $\pMeas(\CE)$ to denote the space of probability
measures on a Polish space $\CE$, endowed with the topology of weak convergence.

\begin{definition}\label{def:Feller}
A {\em transition semigroup} $\CP_t$ on a Polish space $\CE$ is {\em Feller} if it maps  $\CC_b(\CE)$
into  $\CC_b(\CE)$.
\end{definition}

\begin{remark}\label{rem:cont}
This definition is equivalent to the requirement that $x \mapsto \CP_t(x,\cdot\,)$ is continuous from
$\CE$ to $\pMeas(\CE)$. As a consequence, Feller semigroups preserve the weak topology in the
sense that if $\mu_n \to \mu$ in $\pMeas(\CE)$, then $\CP_t \mu_n \to \CP_t \mu$ in $\pMeas(\CE)$ 
for every given $t$.
\end{remark}

Now that we have defined the ``good'' objects for the ``noisy'' part of our construction, we turn
to the trajectories on the state space. We are looking for a space which has good topological
properties but which is large enough to contain most interesting examples. One such space
is the space of {\em \cadlag} paths (continu \`a droite, limite \`a gauche --- continuous on the right,
limits on the left), which can be turned into a Polish space when equipped with a suitable topology.

\begin{definition}\label{def:cadlag}
Given a Polish space $\CE$ and a positive number $T$, the space $\CD([0,T],\CE)$ is the
set of functions $f\from[0,T]\to \CE$ that are right-continuous and whose left-limits exist at every point.
A sequence $\{f_n\}_{n\in\N}$ converges to a limit $f$ if and only if there exists a sequence $\{\lambda_n\}$
of continuous and increasing functions $\lambda_n\from[0,T] \to [0,T]$ satisfying $\lambda_n(0) = 0$,
$\lambda_n(T)=T$, and such that
\begin{equ}[e:cond1]
\lim_{n\to \infty}\sup_{0\le s < t \le T} \Bigl|{\log{\lambda_n(t) - \lambda_n(s) \over t-s}}\Bigr| = 0\;,
\end{equ}
and
\begin{equ}[e:cond2]
\lim_{n\to\infty} \sup_{0\le t\le T} d\sml(f_n(t), f(\lambda_n(t))\smr) = 0\;,
\end{equ}
where $d$ is any totally bounded metric on $\CE$ which generates its topology.

The space $\CD(\R_+, \CE)$ is the space of all functions from $\R_+$ to $\CE$ such that their restrictions to $[0,T]$
are in $\CD([0,T],\CE)$ for all $T>0$. A sequence converges in  $\CD(\R_+, \CE)$ if there exists a sequence
$\{\lambda_n\}$ of continuous and increasing functions $\lambda_n\from\R_+\to \R_+$ satisfying $\lambda_n(0) = 0$
and such that \eref{e:cond1} and \eref{e:cond2} hold.
\end{definition}

It can be shown (see \eg \cite{MR88a:60130} for a proof) that the spaces $\CD([0,T],\CE)$ and $\CD(\R_+, \CE)$ 
are Polish when equipped with the above topology (usually called the Skorohod topology). Notice that the 
space $\CD([0,T],\CE)$ has a natural embedding into $\CD(\R_+, \CE)$ by setting $f(t) = f(T)$ for $t>T$ and
that this embedding is continuous. However, the restriction operator from $\CD(\R_+, \CE)$ to $\CD([0,T],\CE)$
is not continuous, since the topology on $\CD([0,T],\CE)$ imposes that $f_n(T) \to f(T)$, which is not imposed
by the topology on $\CD(\R_+, \CE)$.

In many interesting situations, it is enough to work with continuous sample paths, which live in much simpler spaces:

\begin{definition}\label{def:continuous}
Given a Polish space $\CE$ and a positive number $T$, the space $\CC([0,T],\CE)$ is the
set of continuous functions $f\from[0,T]\to \CE$ equipped with the supremum norm.

The space $\CC(\R_+, \CE)$ is the space of all functions from $\R_+$ to $\CE$ such that their restrictions to $[0,T]$
are in $\CC([0,T],\CE)$ for all $T>0$. A sequence converges in  $\CC(\R_+, \CE)$ if all its restrictions converge.
\end{definition}

It is a standard result that the spaces $\CC([0,T],\CE)$ and $\CC(\R_+, \CE)$ are Polish if $\CE$ is Polish. We 
can now turn to the definition of the systems we are interested in.

\subsection{Definition of a SDS}

Let us recall the following standard notations.
Given a product space $\State \times \Noise$, we denote by $\Pi_\State$ and $\Pi_\Noise$ the maps that select the
first (resp.\ second) component of an element. Also, given two measurable spaces $\CE$ and $\CF$, a measurable
map $f\from\CE\to\CF$, and a measure $\mu$ on $\CE$, we define the measure $f^*\mu$ on $\CF$ in the natural 
way by $f^*\mu = \mu \circ f^{-1}$.

We first define the class of noise processes we will be interested in:
\begin{definition}\label{def:noise}
A quadruple $\sml(\Noise, \{\CP_t\}_{t\ge0}, \pNoise, \{\theta_t\}_{t\ge 0}\smr)$ is called a 
{\em stationary noise process} if it satisfies the following:
\begin{claim}[iii)]
\item[i)] $\Noise$ is a Polish space,
\item[ii)] $\CP_t$ is a Feller transition semigroup on $\Noise$, which accepts $\pNoise$ as its unique
invariant measure,
\item[iii)] The family $\{\theta_t\}_{t>0}$ is a semiflow of measurable maps on $\Noise$ satisfying the property
$\theta_t^*\CP_t^{}(x,\cdot) = \delta_x$ for every $x\in\Noise$.
\end{claim}
\end{definition}

This leads to the following definition of SDS, which is intentionally kept as close as possible to 
the definition of RDS in \cite[Def.~1.1.1]{RDS}:
\begin{definition}\label{def:SDS}
A {\em stochastic dynamical system} on the Polish space $\State$ over the stationary noise
process $\sml(\Noise, \{\CP_t\}_{t\ge0}, \pNoise, \{\theta_t\}_{t\ge 0}\smr)$ is a mapping
\begin{equ}
\phi\from \R_+ \times \State \times \Noise \to \State \;,\qquad (t,x,w)\mapsto \phi_t(x,w)\;,
\end{equ}
with the following properties:
\begin{claim}[(SDS2)]
\item[\maketag{SDS1}{tag:SDS1}] {\em Regularity of paths:} For every $T>0$, $x\in\State$, and $w\in\Noise$, the map $\Phi_T(x,w)\from[0,T]\to \State$ defined by
\begin{equ}{}
\Phi_T(x,w)(t) = \phi_t(x,\theta_{T-t}w)\;,
\end{equ}
belongs to $\CD([0,T],\State)$.
\item[\maketag{SDS2}{tag:SDS2}] {\em Continuous dependence:} The maps $(x,w)\mapsto \Phi_T(x,w)$ are continuous from
$\State\times\Noise$ to $\CD([0,T],\State)$ for every $T>0$.
\item[\maketag{SDS3}{tag:SDS3}] {\em Cocycle property:} The family of mappings $\phi_t$ satisfies
\begin{equs}
\phi_0(x,w) &= x\;,\\
\phi_{s+t}(x,w) &= \phi_s\sml(\phi_t(x,\theta_s w),w\smr)\;, \label{e:cocycle}
\end{equs}
for all $s,t > 0$, all $x\in\State$, and all $w\in\Noise$.
\end{claim}
\end{definition}

\begin{remark}\label{rem:Hans}
The above definition is very close to the definition of {\em Markovian random dynamical system}
introduced in \cite{MR93c:60096}. Beyond the technical differences, the main difference is a shift in
the viewpoint: a Markovian RDS is built on top of a RDS, so one can analyse it from both a semigroup
point of view and a RDS point of view. In the case of a SDS as defined above, there is no underlying RDS
(although one can always construct one), so the semigroup point of view is the only one we consider.
\end{remark}

\begin{remark}
The cocycle property \eref{e:cocycle} looks different from the cocycle property for random dynamical systems.
Actually, in our case $\phi$ is a {\em backward cocycle} for $\theta_t$, which is reasonable since,
as a ``left inverse'' for $\CP_t$, $\theta_t$ actually pushes time backward.
Notice also that, unlike in the definition of RDS, we require some continuity property with respect to 
the noise to hold. This continuity property sounds quite restrictive, but it is actually mainly a matter of
choosing a topology on $\Noise$, which is in a sense ``compatible'' with the topology on $\State$.
\end{remark}

Similarly, we define a continuous (where ``continuous'' should be thought of as continuous with respect
to time) SDS by

\begin{definition}\label{def:contSDS}
A SDS is said to be {\em continuous} if $\CD([0,T],\State)$ can be replaced by $\CC([0,T],\State)$ in
the above definition.
\end{definition}

\begin{remark}\label{rem:compat}
One can check that the embeddings $\CC([0,T],\State) \hookrightarrow \CD([0,T],\State)$
and $\CC(\R_+,\State) \hookrightarrow \CD(\R_+,\State)$ are continuous, so a continuous SDS also satisfies 
Definition~\ref{def:SDS} of a SDS.
\end{remark}

Given a SDS as in Definition~\ref{def:SDS} and an initial condition $x_0 \in \State$, we now turn to the
construction of a stochastic process with initial condition $x_0$ constructed in a natural way from $\phi$. 
First, given $t\ge 0$ and $(x,w) \in \State\times\Noise$,
we construct a probability measure $\CQ_t(x,w;\cdot\,)$ on $\State\times\Noise$ by
\begin{equ}[e:defQt]
\CQ_t(x,w; A\times B) = \int_B \delta_{\phi_t(x,w')}(A)\, \CP_t(w, dw')\;,
\end{equ}
where $\delta_x$ denotes the delta measure located at $x$. The following result is elementary:
\begin{lemma}\label{lem:semigroup}
Let $\phi$ be a SDS on $\State$ over $\sml(\Noise, \{\CP_t\}_{t\ge0}, \pNoise, \{\theta_t\}_{t\ge 0}\smr)$ and define the family
of measures $\CQ_t(x,w;\cdot\,)$ by \eref{e:defQt}. Then $\CQ_t$ is a Feller transition semigroup on $\State \times \Noise$.
Furthermore, it has the property that if $\Pi_\Noise^* \mu = \pNoise$
for a measure $\mu$ on $\State \times\Noise$, then $\Pi_\Noise^* \CQ_t \mu = \pNoise$.
\end{lemma}
\begin{proof}
The fact that  $\Pi_\Noise^* \CQ_t \mu = \pNoise$
follows from the invariance of $\pNoise $ under $\CP_t$.
We now check that $\CQ_t$ is a Feller transition semigroup. Conditions {\it i)} and {\it iii)} follow immediately from 
the properties of $\phi$. The continuity of $\CQ_t(x,w;\cdot\,)$ with respect to $(x,w)$ is a straightforward consequence
of the facts that $\CP_t$ is Feller and that $(x,w) \mapsto \phi_t(x,w)$ is continuous (the latter statement follows from
\eref{tag:SDS2} and the definition of the topology on $\CD([0,t],\State)$).

It thus remains only to check that the Chapman-Kolmogorov equation holds. We have from the cocycle property:
\begin{equs}
\CQ_{s+t}(x,w;&A\times B) =  \int_B \delta_{\phi_{s+t}(x,w')}(A)\, \CP_{s+t}(w, dw') \\
&=  \int_B \int_\State \delta_{\phi_{s}(y,w')}(A) \delta_{\phi_{t}(x,\theta_s w')}(dy)\, \CP_{s+t}(w, dw')\\
&=  \int_\Noise \int_B \int_\State \delta_{\phi_{s}(y,w')}(A) \delta_{\phi_{t}(x,\theta_s w')}(dy)\, \CP_{s}(w'', dw')\, \CP_{t}(w, dw'')\;.
\end{equs}
The claim then follows from the property $\theta_s^*\CP_s(w'',dw') = \delta_{w''}(dw')$ by exchanging the order of integration.
\end{proof}

\begin{remark}\label{rem:manypoint}
Actually, \eref{e:defQt} defines the evolution of the one-point process generated by $\phi$. The $n$-points process
would evolve according to
\begin{equ}
\CQ_t^{(n)}(x_1,\ldots,x_n,w; A_1\times\ldots\times A_n \times B) = \int_B \prod_{i=1}^n \delta_{\phi_t(x_i,w')}(A_i )\, \CP_t(w, dw')\;.
\end{equ}
One can check as above that this defines a Feller transition semigroup on $\State^n \times \Noise$.
\end{remark}

This lemma suggests the following definition:
\begin{definition}\label{def:IC}
Let $\phi$ be a SDS as above. Then a probability measure $\mu$ on $\State\times\Noise$ is called a {\em generalised 
initial condition} for $\phi$ if $\Pi_\Noise^* \mu = \pNoise$. We denote by $\Meas$ the space of generalised initial conditions
endowed with the topology of weak convergence. Elements of $\Meas$ that are of the form
$\mu = \delta_x\times \pNoise$ for some $x\in\State$ will be called {\em initial conditions}.
\end{definition}

Given a generalised initial condition $\mu$, it is natural construct a stochastic process $(x_t,w_t)$ on $\State\times\Noise$ by
drawing its initial condition according to $\mu$ and then evolving it according to the transition semigroup $\CQ_t$. 
The marginal $x_t$ of this process on $\State$ will be called the {\em process generated by $\phi$ for $\mu$}. We will denote
by $\Evol \mu$ the law of this process (\ie $\Evol\mu$ is a measure on $\CD(\R_+,\State)$ in the general case and
a measure on $\CC(\R_+,\State)$ in the continuous case). More rigorously, we define for every $T>0$ the measure
$\Evol_T\mu$ on $\CD([0,T],\State)$ by
\begin{equ}
\Evol_T\mu = \Phi_T^*\CP_t \mu\;,
\end{equ}
where $\Phi_T$ is defined as in \eref{tag:SDS1}. By the embedding  $\CD([0,T],\State)\hookrightarrow\CD(\R_+,\State)$,
this actually gives a family of measures on $\CD(\R_+,\State)$.
It follows from the cocycle property that the restriction to $\CD([0,T],\State)$
of $\Evol_{T'}\mu$ with $T'>T$ is equal to $\Evol_T\mu$. The definition of the topology on $\CD(\R_+,\State)$ 
does therefore imply that the sequence $\Evol_T\mu$ converges weakly to a unique measure on $\CD(\R_+,\State)$
that we denote by $\Evol\mu$. A similar argument, combined with \eref{tag:SDS2} yields
\begin{lemma}\label{lem:cont}
Let $\phi$ be a SDS. Then, the operator $\Evol$ as defined above is continuous from $\Meas$ to 
$\pMeas(\CD(\R_+,\State))$. \qed
\end{lemma}

This in turn motivates the following equivalence relation:
\begin{definition}\label{def:equiv}
Two generalised initial conditions $\mu$ and $\nu$ of a SDS $\phi$ are {\em equivalent} if the processes generated
by $\mu$ and $\nu$ are equal in law. In short, $\mu \sim \nu \Leftrightarrow \Evol \mu = \Evol \nu$.
\end{definition}
The physical interpretation of this notion of equivalence is that the noise space contains some redundant information
that is not required to construct the future of the system. Note that 
this does not necessarily mean that the noise space could be reduced in order to have a more ``optimal'' description
of the system. For example, if the process $x_t$ generated by any generalised initial condition is Markov, then all the information
contained in $\Noise$ is redundant in the above sense (\ie $\mu$ and $\nu$ are equivalent if $\Pi_\State^*\mu = \Pi_\State^*\nu$).
This does of course not mean that $\Noise$ can be entirely thrown away in the above description (otherwise, 
since the map $\phi$ is deterministic, the evolution
would become deterministic).

The main reason for introducing the notion of SDS is to have a framework in which one can study 
ergodic properties of physical systems with memory. It should be noted that it is designed to describe 
systems where the memory is \textit{extrinsic}, as opposed to systems with \textit{intrinsic} memory like
stochastic delay equations.
We present in the next subsection a few elementary ergodic results in
the framework of SDS.

\subsection{Ergodic properties}

In the theory of Markov processes, the main tool for investigating ergodic properties is the {\em invariant measure}.
In the setup of SDS, we say that a measure $\mu$ on $\State\times\Noise$ is invariant for the SDS $\phi$ if it
is invariant for the Markov transition semigroup $\CQ_t$ generated by $\phi$. We say that a 
measure $\mu$ on $\State\times\Noise$ is {\em stationary} for $\phi$ if one has
\begin{equ}
 \CQ_t \mu \sim  \mu\;,\quad \forall t>0\;,
\end{equ}
\ie if the process on $\State$ generated by $\mu$ is stationary. Following our philosophy of considering only what happens
on the state space $\State$, we should be interested in stationary measures, disregarding completely whether they
are actually invariant or not. In doing so, we could be afraid of loosing many convenient results from the
well-developed
theory of Markov processes. Fortunately, the following lemma shows that the set of invariant 
measures and the set of stationary measures are actually the same, when quotiented 
by the equivalence relation of Definition~\ref{def:equiv}.
\begin{proposition}\label{prop:invariant}
Let $\phi$ be a SDS and let $\mu$ be a stationary measure for $\phi$. Then, there exists
a measure $\mu_\star \sim \mu$ which is invariant for $\phi$.
\end{proposition}

\begin{proof}
Define the ergodic averages
\begin{equ}[e:ergodic]
\CR_T \mu = {1 \over T}\int_0^T \CQ_t \mu\, dt\;.
\end{equ}
Since $\mu$ is stationary, we have $\Pi_\State^*\CR_T \mu = \Pi_\State^* \mu$ for every $T$.
Furthermore, $\Pi_\Noise^*\CR_T \mu = \pNoise$ for every $T$, therefore the sequence of measures
$\CR_T \mu$ is tight on $\State \times \Noise$. Let $\mu_\star$ be any of its accumulation points
in $\pMeas(\State \times \Noise)$. Since  $ \CQ_t $ is Feller, $\mu_\star$ is invariant for $ \CQ_t $
and,  by \lem{lem:cont}, one has $\mu_\star \sim \mu$.
\end{proof}

From a mathematical point of view, it may in some cases be interesting  to know whether the invariant
measure $\mu_\star$ constructed in \prop{prop:invariant} is uniquely determined by $\mu$. From
an intuitive point of view, this uniqueness property should hold
if the information contained in the trajectories on the state space $\State $ is sufficient to 
reconstruct the evolution of the noise. This intuition is made rigorous by the following proposition.

\begin{proposition}\label{prop:faithful}
Let $\phi$ be a SDS, define ${\cscr W}_T^x$ as the $\sigma$-field on $\Noise$ generated
by the map $\Phi_T(x,\cdot\,)\from\Noise\to\CD([0,T],\State)$, and set ${\cscr W}_T = \bigwedge_{x\in\State} {\cscr W}_T^x$. 
Assume that ${\cscr W}_T \subset {\cscr W}_{T'}$ for $T < T'$ and that
${\cscr W} = \bigvee_{T\ge 0} {\cscr W}_T$ is equal to the Borel 
$\sigma$-field on $\Noise$. Then, for $\mu_1$ and $\mu_2$ two invariant measures, one has
the implication $\mu_1 \sim \mu_2 \Rightarrow \mu_1 = \mu_2$.
\end{proposition}

\begin{proof}
Assume $\mu_1 \sim \mu_2$ are two invariant measures for $\phi$. 
Since ${\cscr W}_T \subset {\cscr W}_{T'}$ if $T < T'$, their equality follows if one can show that, for every $T > 0$,
\begin{equ}[e:equal]
\expect \bigl(\mu_1\,|\, {\cscr X} \otimes {\cscr W}_T \bigr) = \expect \bigl(\mu_2\,|\, {\cscr X} \otimes {\cscr W}_T\bigr)\;,
\end{equ}
where ${\cscr X}$ denotes the Borel $\sigma$-field on $\State$. 

Since $\mu_1 \sim \mu_2$, one has in particular $\Pi_\State^*\mu_1 = \Pi_\State^*\mu_2$, so let us call 
this measure $\nu$. Since $\Noise$ is Polish, we then have the disintegration $x \mapsto \mu_i^x$, yielding formally
$\mu_i(dx,dw) = \mu_i^x(dw)\,\nu(dx)$, where $\mu_i^x$ are probability measures on $\Noise$. 
(See \cite[p.~196]{MR58:18632} for a proof.) Fix $T>0$ and define the family $\mu_i^{x,T}$ of probability measures
on $\Noise$ by
\begin{equ}
\mu_i^{x,T} = \int_\Noise \CP_t(w,\cdot\,)\,\mu_i^x(dw)\;.
\end{equ}
With this definition, one has
\begin{equ}
\Evol_T\mu_i = \int_\State \bigl(\Phi_T(x,\cdot\,)^*\mu_i^{x,T}\bigr)\,\nu(dx)\;.
\end{equ}
Let $e_0\from \CD([0,T],\State)\to\State$ be the evaluation map at $0$, then
\begin{equ}
\expect \bigl(\Evol_T\mu_i\,|\, e_0 = x\bigr) =  \bigl(\Phi_T(x,\cdot\,)^*\mu_i^{x,T}\bigr)\;,
\end{equ}
for $\nu$-almost every $x\in\State$.
Since $\Evol_T\mu_1 = \Evol_T\mu_2$, one therefore has
\begin{equ}[e:equal_2]
\expect \bigl(\mu_1^{x,T} \,|\, {\cscr W}_T^x \bigr) 
= \expect \bigl(\mu_2^{x,T} \,|\, {\cscr W}_T^x \bigr)\;,
\end{equ}
for  $\nu$-almost every $x\in\State$. On the other hand, the invariance of $\mu_i$ implies that,
 for every $A\in {\cscr X}$ and every $B \in {\cscr W}_T$, one has the equality
\begin{equ}
\mu_i(A \times B) = \int_\State \int_B  \chi_A \bigl(\phi_T(x,w)\bigr)\, \mu_i^{x,T}(dw)\,\nu(dx)\;.
\end{equ}
Since $\phi_T(x,\cdot\,)$ is ${\cscr W}_T^x$-measurable and $B \in {\cscr W}_T^x$, this is equal to
\begin{equ}
\int_\State \int_B  \chi_A \bigl(\phi_T(x,w)\bigr)\, \expect \bigl(\mu_i^{x,T}\,|\, {\cscr W}_T^x\bigr)(dw)\,\nu(dx)\;.
\end{equ}
Thus \eref{e:equal_2} implies \eref{e:equal} and the proof of \prop{prop:faithful} is complete.
\end{proof}

The existence of an invariant measure is usually established by finding a Lyapunov function.
In this setting, Lyapunov functions are given by the following definition.

\begin{definition}\label{def:Lyap}
Let $\phi$ be a SDS and let $F\from\State \to [0,\infty)$ be a continuous function. Then $F$
is a {\em Lyapunov function for $\phi$} if it satisfies the following conditions:
\begin{claim}[(L2)]
\item[(L1)] The set $F^{-1}([0,C])$ is compact for every $C \in [0,\infty)$.
\item[(L2)] There exist constants $C$ and $\gamma > 0$ such that
\begin{equ}[e:condLyap]
\int_{\State\times\Noise} F(x)\,\sml(\CQ_t \mu\smr)(dx,dw)  \le C + e^{-\gamma t} \int_{\State} F(x)\,\sml(\Pi_\State^*\mu\smr)(dx)\;,
\end{equ}
for every $t>0$ and every generalised initial condition $\mu$ such that the right-hand side is finite.
\end{claim}
\end{definition}

It is important to notice that one does {\em not} require $F$ to be a Lyapunov function for the
transition semigroup $\CQ_t$, since \eref{e:condLyap} is only required to hold for measures
$\mu$ satisfying $\Pi_\Noise^*\mu = \pNoise$. One nevertheless has the following result:

\begin{lemma}\label{lem:Lyap}
Let $\phi$ be a SDS. If there exists a Lyapunov function $F$ for $\phi$, then there exists also
an invariant measure $\mu_\star$ for $\phi$, which satisfies
\begin{equ}[e:estInv]
\int_{\State\times\Noise} F(x)\,\mu_\star(dx,dw) \le C\;.
\end{equ} 
\end{lemma}

\begin{proof}
Let $x\in\State$ be an arbitrary initial condition, set $\mu = \delta_x \times \pNoise$, and
define the ergodic averages $\CR_T \mu$ as in \eref{e:ergodic}. Combining (L1) and (L2) with the fact that
$\Pi_\Noise^*\CR_T \mu = \pNoise$, one immediately gets the tightness of the sequence $\{\CR_T \mu\}$. By the 
standard Krylov-Bogoloubov argument, any limiting point of $\{\CR_T \mu\}$ is an invariant measure for $\phi$.
The estimate \eref{e:estInv} follows from \eref{e:condLyap}, combined with the fact that $F$ is continuous.
\end{proof}

This concludes our presentation of the abstract framework in which we analyse the
ergodic properties of \eref{e:mainequ}.

\section{Construction of the SDS}
\label{sec:setting}

In this section, we construct a continuous stochastic dynamical system which
yields the solutions to \eref{e:mainequ} in an appropriate sense.

First of all, let us discuss what we mean by ``solution'' to \eref{e:mainequ}. 
\begin{definition}\label{def:solution}
Let $\{x_t\}_{t\ge 0}$ be a stochastic process with continuous sample paths. We say
that $x_t$ is a {\em solution} to \eref{e:mainequ} if the stochastic process $N(t)$ defined by
\begin{equ}[e:fBmsol]
N(t) = x_t - x_0 - \int_0^t f(x_s)\,ds\;,
\end{equ}
is equal in law to $\sigma B_H(t)$, where $\sigma$ is as in \eref{e:mainequ} and $B_H(t)$ is
a $n$-dimensional fBm with Hurst parameter $H$.
\end{definition}

We will set up our SDS in such a way that, for every generalised initial condition $\mu$, the canonical
process associated to the measure $\Evol\mu$ is a solution to \eref{e:mainequ}. This will be the content
of \prop{prop:defSDS} below. In order to achieve this,
our main task is to set up a noise process in a way which complies to Definition~\ref{def:noise}.

\subsection{Representation of the fBm}

In this section, we give a representation of the fBm $B_H(t)$ with 
Hurst parameter $H \in (0,1)$ which is suitable for our analysis. 
Recall that, by definition, $B_H(t)$ is a centred Gaussian process satisfying $B_H(0) = 0$ and
\begin{equ}[e:propfBm]
 \expect|B_H(t) - B_H(s)|^2 = |t-s|^{2H}\;.
\end{equ}
Naturally, a {\it two-sided fractional 
Brownian motion} by requiring that \eref{e:propfBm} holds for all $s,t\in\R$. 
Notice that, unlike for the normal Brownian motion, the two-sided fBm is \textit{not} obtained by gluing 
two independent
copies of the one-sided fBm together at $t=0$.
We have the following useful representation of the two-sided fBm, which is also (up to the normalisation constant)
the representation used in the original paper \cite{MR39:3572}.

\begin{lemma}\label{lem:fbm}
Let $w(t)$, $t\in\R$ be a two-sided Wiener process and let $H \in (0,1)$. Define for some constant $\alpha_H$ the process
\begin{equ}[e:repres]
B_H(t) =  \alpha_H\int_{-\infty}^0 \Kern(-r)\,\sml(dw(r+t) - dw(r)\smr)\;.
\end{equ}
Then there exists a choice of $\alpha_H$ such that $B_H(t)$ is a two-sided fractional Brownian motion with 
Hurst parameter $H$. \qed
\end{lemma}

\begin{notation}
Given the representation \eref{e:repres} of the fBm with Hurst parameter $H$, we call $w$ the ``Wiener process associated to $B_H$''. We
also refer to $\{w(t)\,:\, t\le 0\}$ as the ``past'' of $w$ and to $\{w(t)\,:\, t > 0\}$ as the ``future'' of $w$. We similarly refer to the ``past'' and the ``future'' of $B_H$. Notice the notion of future for $B_H$ is different from the notion of future for $w$ in terms of $\sigma$-algebras, since the future of $B_H$
depends on the past of $w$.
\end{notation}

\begin{remark}
The expression \eref{e:repres} looks strange at first sight, but one should actually think of $B_H(t)$ as being given by 
$B_H(t) = \tilde B_H(t) - \tilde B_H(0)$, where
\begin{equ}[e:Btilde]
\raise0.4em\hbox{``}\tilde B_H(t) = \alpha_H \int_{-\infty}^t \Kern(t-s)\,dw(s)\;.\;\raise0.4em\hbox{''}
\end{equ}
This expression is strongly reminiscent of the usual representation of the stationary Ornstein-Uhlenbeck process, but with an algebraic kernel instead of an exponential one.
Of course, \eref{e:Btilde} does not make any sense since $(t-s)^{H-{1\over 2}}$ is not square integrable. Nevertheless, \eref{e:Btilde} 
has the advantage of explicitly showing the stationarity of the increments for the two-sided fBm. 
\end{remark}


\subsection{Noise spaces}

In this section, we introduce the family of spaces that will be used to model our noise.
Denote by $\CC_0^\infty(\R_-)$ the set 
of $\CC^\infty$ function $w\from(-\infty,0] \to \R$ satisfying $w(0) = 0$ and having compact support. 
Given a parameter $H\in(0,1)$, 
we define for every $w\in\CC_0^\infty(\R_-)$ the norm
\begin{equ}[e:defN]
\|w\|_H = \sup_{t,s\in\R_-} {|w(t)-w(s)| \over |t-s|^{1-H\over 2} \sml(1+|t| + |s|\smr)^{1\over2}}\;.
\end{equ}
We then define the Banach space $\CH_H$ to be the closure of $\CC_0^\infty(\R_-)$ under the norm
$\|\nb\cdot\nb\|_H$. The following lemma is important in view of the framework exposed in Section~\ref{sec:general}:

\begin{lemma}\label{lem:separ}
The spaces $\CH_H$ are separable.
\end{lemma}

\begin{proof}
It suffices to find a norm $\|\cdot\|_\star$ which is stronger than $\|\cdot\|_H$ and such that
the closure of $\CC_0^\infty(\R_-)$ under $\|\cdot\|_\star$ is separable. One example of such a norm
is given by $\|w\|_\star= \sup_{t < 0} |t \dot w(t)|$.
\end{proof}

Notice that it is crucial to define $\CH_H$ as the closure of $\CC_0^\infty$ under $\|\cdot\|_H$. If we defined
it simply as the space of all functions with finite $\|\cdot\|_H$-norm, it would not be separable. (Think of the 
space of bounded continuous functions, versus the space of continuous functions vanishing at infinity.)

In view of the representation \eref{e:repres}, we define the linear operator $\CD_H$ on functions $w\in\CC_0^\infty$ by
\begin{equ}[e:defD]
\sml(\CD_H w\smr)(t) =  \alpha_H\int_{-\infty}^0 (-s)^{H-{1\over 2}}\sml(\dot w(s+t)- \dot w(s)\smr)\,ds\;,
\end{equ}
where $\alpha_H$ is as in \lem{lem:fbm}.
We have the following result:

\begin{lemma}\label{lem:spaces}
Let $H \in (0,1)$ and let $\CH_H$ be as above. Then the operator $\CD_H$, formally defined by \eref{e:defD},
 is continuous from
$\CH_H$ into $\CH_{1-H}$. Furthermore, the operator $\CD_H$ has a bounded inverse, given by the formula
\begin{equ}
\CD_H^{-1} = \gamma_H \CD_{1-H}\;,
\end{equ}
for some constant $\gamma_H$ satisfying $\gamma_H = \gamma_{1-H}$.
\end{lemma}

\begin{remark}
The operator $\CD_H$ is actually (up to a multiplicative constant) a fractional integral of order $H-{1\over 2}$ which
 is renormalised in such a way that one gets rid of the divergence at $-\infty$. It is therefore not surprising that 
the inverse of $\CD_H$ is $\CD_{1-H}$. 
\end{remark}

\begin{proof}
For $H = {1\over 2}$, $\CD_H$ is the identity and there is nothing to prove. We therefore assume in the sequel that
$H \neq {1\over 2}$. 

We first show that $\CD_H$ is continuous from $\CH_H$ into $\CH_{1-H}$.
One can easily check that $\CD_H$ maps $\CC_0^\infty$ into the set of $\CC^\infty$ functions which
converge to a constant at $-\infty$. This set can be seen to belong to $\CH_{1-H}$ by a 
simple cutoff argument, so it suffices
to show that $\|\CD_H w\|_{1-H} \le C\|w\|_H$ for $w\in\CC_0^\infty$.
Assume without loss of generality that $t>s$ and define $h = t-s$. We then have
\begin{equs}
\sml(\CD_H w\smr)(t)- \sml(\CD_H w\smr)(s) &= \alpha_H \int_{-\infty}^{s} \bigl(\Kern(t-r) - \Kern(s-r)\bigr)\,dw(r) \\
	&\quad + \alpha_H\int_s^t\Kern(t-r)\,dw(r) \;.
\end{equs}
Splitting the integral and integrating by parts yields
\begin{equs}
\sml(\CD_H w\smr)(t) &- \sml(\CD_H w\smr)(s) = -\alpha_H\sml(H-{\textstyle{1\over 2}}\smr) \int_{s-h}^s (s-r)^{H-{3\over 2}}\sml(w(r)-w(s)\smr)\,dr \\
& +\alpha_H\sml(H-{\textstyle{1\over 2}}\smr) \int_{t-2h}^t (t-r)^{H-{3\over 2}}\sml(w(r)-w(t)\smr)\,dr \\
& +\alpha_H\sml(H-{\textstyle{1\over 2}}\smr)\int_{-\infty}^{s-h} \bigl((t-r)^{H-{3\over 2}} - (s-r)^{H-{3\over 2}}\bigr)\sml(w(r)-w(s)\smr)\,dr \\
& +\alpha_H (2h)^{H-{1\over 2}}\sml(w(t) - w(s)\smr)\\[0.5em]
& \equiv T_1 + T_2 + T_3 + T_4\;.
\end{equs}
We estimate each of these terms separately. For $T_1$, we have
\begin{equ}
|T_1| \le C \sml(1 + |s| + |t|\smr)^{1/2} \int_0^h r^{H-{3\over 2}+ {1-H\over 2}}\,dr \le C h^{H\over2} \sml(1 + |s| + |t|\smr)^{1/2}\;.
\end{equ}
The term $T_2$ is bounded by $C h^{H\over2} \sml(1 + |s| + |t|\smr)^{1/2}$ in a similar way. Concerning $T_3$, we bound it by
\begin{equs}
|T_3| &\le C \int_h^\infty \bigl(r^{H-{3\over 2}} - (h+r)^{H-{3\over 2}}\bigr)\sml(w(s-r) - w(s)\smr)\,dr \\
 &\le C h \int_h^\infty r^{H-{5\over 2}}r^{1-H\over 2} \sml(1 + |s| + |r|\smr)^{1/2}\,dr \\
& \le C h^{H \over 2} \sml(1 + |s|\smr)^{1/2} + C h \int_h^\infty r^{{H\over 2}-2} \sml(h+r\smr)^{1/2}\,dr\\
& \le C h^{H \over 2} \sml(1 + |s| + h \smr)^{1/2}  \le C h^{H \over 2} \sml(1 + |s| + |t|\smr)^{1/2}\;.
\end{equs}
The term $T_4$ is easily bounded by $C h^{H \over 2} \sml(1 + |s| + |t|\smr)^{1/2}$, using the fact that  $w \in \CH_H$. 
This shows that $\CD_H$ is bounded from $\CH_H$ to $\CH_{1-H}$. 

It remains to show that $\CD_H \circ \CD_{1-H}$
is a multiple of the identity. For this, notice that if $w \in \CC_0^\infty$, then one has in the
notations of \cite[pp.~94--95]{Frac} the following identities
\begin{equs}[2]
\bigl(\CD_H w\bigr)(t) &= -\alpha_H \Gamma(H + {\textstyle{1\over 2}}) \bigl(\bigl(I_+^{H-{1\over 2}} w\bigr)(t) - \bigl(I_+^{H-{1\over 2}} w\bigr)(0)\bigr)\;,&\quad H &> {\textstyle{1\over 2}}\;,\\
\bigl(\CD_H w\bigr)(t) &= -\alpha_H \Gamma(H + {\textstyle{1\over 2}}) \bigl(\bigl(D_+^{{1\over 2}-H} w\bigr)(t) - \bigl(D_+^{{1\over 2}-H} w\bigr)(0)\bigr)\;,&\quad H &< {\textstyle{1\over 2}}\;.
\end{equs}
Furthermore, \eref{e:defD} shows that $\CD_H w = 0$ if $w$ is a constant. The claim then follows immediately
from the fact that if $w \in \CC_0^\infty$ and $\alpha \in (0,1)$, one has $D_+^\alpha I_+^\alpha w  = w$ and $I_+^\alpha  D_+^\alpha w= w$ (see \cite[Thm.~2.4]{Frac}).
\end{proof}

Since we want to use the operators $\CD_H$ and $\CD_{1-H}$ to switch between Wiener processes and fractional
Brownian motions, it is crucial to show that the sample paths of the two-sided Wiener process belong to every 
$\CH_H$ with probability $1$. Actually, what we show is that the Wiener measure can be constructed as
a Borel measure on $\CH_H$.

\begin{lemma}\label{lem:wienerspace}
There exists a unique Gaussian measure $\Wien$ on $\CH_H$ which is such that the canonical process associated
to it is a time-reversed Brownian motion.
\end{lemma}

\begin{proof}
We start by showing that the $\CH_H$-norm of the Wiener paths has bounded moments of all orders.
It follows from a generalisation of the Kolmogorov criterion \cite[Theorem 2.1]{Yor} that
\begin{equ}[e:locW]
\expect \biggl(\sup_{s,t \in [0,2]} {|w(s) - w(t)| \over |s-t|^{{1-H}\over 2}}\biggr)^p < \infty
\end{equ}
for all $p>0$.
Since the increments of $w$ are independent, this implies that, for every $\eps > 0$, there exists
a random variable $C_1$ such that
\begin{equ}[e:proplocw]
\sup_{|s-t| \le 1} {|w(s) - w(t)| \over |s-t|^{{1-H}\over 2} \sml(1 + |t| + |s|\smr)^{\eps}} < C_1\;,
\end{equ}
with probability $1$, and that all the moments of $C_1$ are bounded. We can therefore 
safely assume in the sequel that $|t-s| > 1$. It follows 
immediately from \eref{e:proplocw} and the triangle inequality that there exists a constant $C$ such that
\begin{equ}[e:proplocw2]
|w(s) - w(t)| \le C C_1 |t-s| \sml(1 + |t| + |s|\smr)^{\eps}\;,
\end{equ}
whenever $|t-s| > 1$. Furthermore, it follows from the time-inversion property of the Brownian motion, 
combined with \eref{e:locW},
that $|w|$ does not grow much faster than $|t|^{1/2}$
for large values of $t$. In particular, for every $\eps' >0$, there exists a random variable $C_2$ such that
\begin{equ}[e:propglobw]
|w(t)| \le C_2\sml(1+|t|\smr)^{{1 \over 2} + \eps'}\;,\qquad \forall t\in\R\;,
\end{equ}
and that all the moments of $C_2$ are bounded.
Combining \eref{e:proplocw2} and \eref{e:propglobw}, we get (for some other constant $C$)
\begin{equ}
|w(s) - w(t)| \le C C_1^{1-H\over 2} C_2^{1+H \over 2} |t-s|^{1-H\over 2} \sml(1+|s| + |t|\smr)^{{H+1 \over 4} + \eps{1-H\over 2} + \eps'{1+H\over 2}}\;.
\end{equ}
The claim follows by choosing for example $\eps = \eps' = (1-H)/4$.

This is not quite enough, since we want the sample paths to belong to the closure of $\CC_0^\infty$ 
under the norm
$\|\cdot\|_H$. Define the function
\begin{equ}
(s,t) \mapsto \Gamma(s,t) = {\sml(1+|t| + |s| \smr)^{2} \over |t-s|} \;.
\end{equ}
By looking at the above proof, we see that we actually proved the stronger statement that 
for every $H\in(0,1)$, one can find a $\gamma > 0$ such that
\begin{equ}
\|w\|_{H,\gamma} = \sup_{s,t}  {\Gamma(s,t)^\gamma |w(s) - w(t)| \over |s-t|^{{1-H}\over 2} \sml(1 + |t| + |s|\smr)^{1\over 2}} < \infty
\end{equ}
with probability $1$. 
Let us call $\CH_{H,\gamma}$ the 
Banach space of functions with finite $\|\cdot\|_{H,\gamma}$-norm. We will show that one has the continuous inclusions:
\begin{equ}[e:inclSpaces]
\CH_{H,\gamma} \hookrightarrow \CH_H \hookrightarrow \CC(\R_-,\R)\;.
\end{equ}
Let us call $\tilde\Wien$ the usual time-reversed Wiener measure on $\CC(\R_-,\R)$ equipped
with the $\sigma$-field ${\cscr R}$ generated by the evaluation functions.
Since $\CH_{H,\gamma}$ is a measurable subset of $\CC(\R_-,\R)$ and $\tilde\Wien(\CH_{H,\gamma})=1$, 
we can restrict $\tilde\Wien$ to a measure on $\CH_H$, equipped with the restriction $\tilde {\cscr R}$
of ${\cscr R}$. It remains to show that $\tilde {\cscr R}$ is equal to the Borel $\sigma$-field ${\cscr B}$ on $\CH_H$.
This follows from the fact that the evaluation functions are ${\cscr B}$-measurable (since they are actually continuous)
and that a countable number of function evaluations suffices to determine the $\|\cdot\|_{H}$-norm of a function.
The proof of \lem{lem:wienerspace} is thus complete if we show \eref{e:inclSpaces}.

Notice first that the function $\Gamma(s,t)$ becomes large
when $|t-s|$ is small or when either $|t|$ or $|s|$ are large, more precisely
we have
\begin{equ}[e:boundGamma]
\Gamma(s,t) > \max\{|s|,\,|t|,\,|t-s|^{-1}\}\;.
\end{equ}
Therefore, functions $w \in \CH_{H,\gamma}$ 
are actually more regular and have better growth properties than what is needed to
have finite $\|\cdot\|_H$-norm. Given $w$ with $\|w\|_{H,\gamma} < \infty$ and 
any $\eps > 0$, we will construct a function
$\tilde w \in \CC_0^\infty$ such that $\|w - \tilde w\|_H < \eps$. Take two $\CC^\infty$ functions
$\phi_1$ and $\phi_2$ with the following shape:
\begin{equ}
\mhpastefig[3/4]{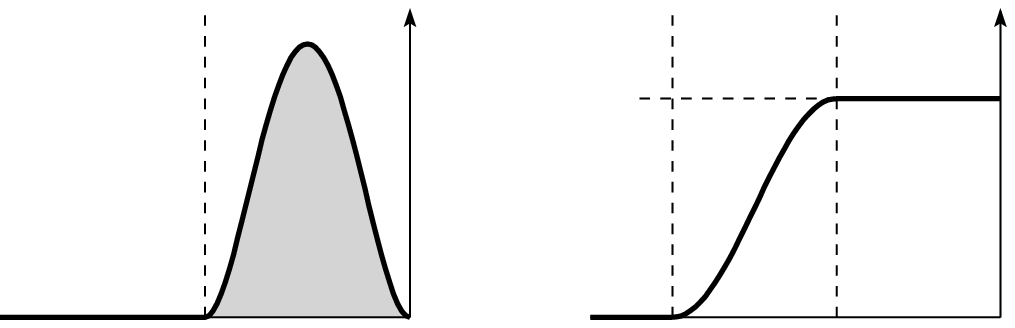}
\end{equ}
Furthermore, we choose them such that:
\begin{equ}
\int_{\R_-} \phi_1(s)\,ds = 1\;,\qquad \Bigl|{d\phi_2(t) \over dt}\Bigr| \le 2\;.
\end{equ}
For two positive constants
$r<1$ and $R>1$ to be chosen later, we define
\begin{equ}
\tilde w(t) = \phi_2(t/R)\int_{\R_-} w(t+s){\phi_1(s/r)\over r}\,ds\;.
\end{equ}
\ie we smoothen out $w$ at length scales smaller than $r$ and we cut it off at 
distances bigger than $R$. A straightforward estimate shows that there exists a constant $C$ such
that
\begin{equ}
\|\tilde w\|_{H,\gamma} \le C\|w\|_{H,\gamma}\;,
\end{equ}
independently of $r<1/4$ and $R>1$.
For $\delta > 0$ to be chosen later, we then divide the
quadrant $K = \{(t,s)\,|\,t,s < 0\}$ into three regions:

\noindent\hspace{1.5cm}
\mhpastefig[4/5]{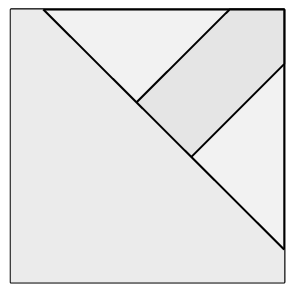}
\begin{minipage}[b]{8cm}
\begin{equs}
K_1 &= \{(t,s)\,|\, |t| + |s| \ge R\} \cap K\;,\\
K_2 &= \{(t,s)\,|\, |t-s| \le \delta\} \cap K \setminus K_1\;,\\
K_3 &= K \setminus (K_1 \cup K_2)\;.
\end{equs}
\end{minipage}
We then bound $\|w-\tilde w\|_H$ by
\begin{equs}
\|w - \tilde w\|_H &\le \sup_{(s,t) \in K_1\cup K_2}{C\|w\|_{H,\gamma} \over \Gamma(t,s)^\gamma} 
+ \sup_{(s,t)\in K_3} {|w(s)-\tilde w(s)| + |w(t)-\tilde w(t)| \over |t-s|^{1-H\over 2}(1+|t|+|s|)^{1\over 2}}\\
&\le C\sml(\delta^\gamma + R^{-\gamma}\smr) \|w\|_{H,\gamma} + 2\delta^{H-1\over 2} \sup_{0<t<R} |w(t)-\tilde w(t)|\;.
\end{equs}
By choosing $\delta$ small enough and $R$ large enough, the first term can be made arbitrarily small. One can then
choose $r$ small enough to make the second term arbitrarily small as well. 
This shows that \eref{e:inclSpaces} holds and therefore the proof of \lem{lem:wienerspace} is complete.
\end{proof}

\subsection{Definition of the SDS}

The results shown so far in this section are sufficient to construct the required SDS. We start by considering the pathwise
solutions to \eref{e:mainequ}. Given a time $T>0$, an initial condition $x\in \R^n$, and a noise $b\in \CC_0([0,T],\R^n)$, we
look for a function $\Phi_T(x,b) \in \CC([0,T],\R^n)$ satisfying
\begin{equ}[e:defSol]
\Phi_T(x,b)(t) = \sigma b(t) + x + \int_0^t f \sml(\Phi_T(x,b)(s)\smr)\,ds\;.
\end{equ}
We have the following standard result:
\begin{lemma}\label{lem:existSol}
Let $f\from\R^n\to\R^n$ satisfy assumptions \textbf{A1} and \textbf{A2}. Then, there exists a unique map
$\Phi_T\from \R^n\times \CC([0,T],\R^n) \to \CC([0,T],\R^n)$ satisfying \eref{e:defSol}. Furthermore, $\Phi_T$
is locally Lipschitz continuous.
\end{lemma}
\begin{proof}
The local (\ie small $T$) existence and uniqueness of continuous solutions to \eref{e:defSol} follows from a standard contraction
argument. In order to show the global existence and the local Lipschitz property, fix $x$, $b$, an $T$, and define
$y(t) = x + \sigma b(t)$. Define $z(t)$ as the solution to the differential equation 
\begin{equ}[e:defz]
\dot z(t) = f(z(t) + y(t))\;,\quad z(0) = 0\;.
\end{equ}
Writing down
the differential equation satisfied by $\|z(t)\|^2$ and using  \textbf{A1} and \textbf{A2}, one sees that \eref{e:defz} possesses
a (unique) solution up to time $T$. One can then set $\Phi_T(x,b)(t) = z(t) + y(t)$ and check that it satisfies \eref{e:defSol}.
The local Lipschitz property of $\Phi_T$ then immediately follows from the local Lipschitz property of $f$.
\end{proof}

We now define the stationary
noise process. For this, we define
$\theta_t\from\CH_H\to\CH_H$ by
\begin{equ}
\sml(\theta_t w\smr)(s) = w(s-t)-w(-t)\;.
\end{equ}
In order to construct the transition semigroup $\CP_t$, we define first $\tilde \CH_H$ like $\CH_H$, but with 
arguments in $\R_+$ instead of $\R_-$, and we write $\tilde \Wien$ for the Wiener measure on $\tilde \CH_H$,
as constructed in \lem{lem:wienerspace} above. Define the function $P_t\from \CH_H \times \tilde \CH_H \to \CH_H$
by
\begin{equ}[e:defTrans]
\sml(P_t(w,\tilde w)\smr)(s) = \cases{\tilde w(t+s) - \tilde w(t) & for $s > -t$,\cr 
	w(t+s) - \tilde w(t) & for $s \le -t$,}
\end{equ}
and set $\CP_t(w,\cdot\,) = P_t(w,\cdot\,)^*\tilde \Wien$. This construction can be visualised by the following picture:
\begin{equ}
\mhpastefig[2/3]{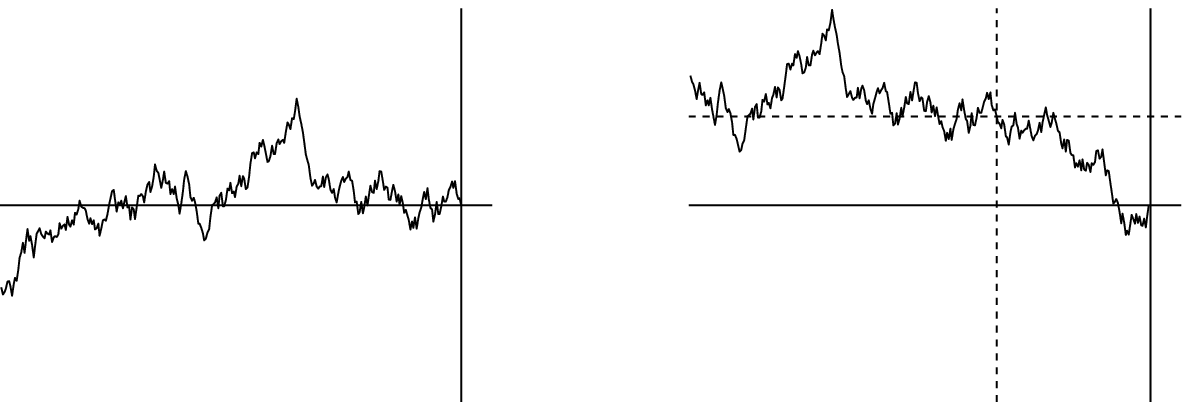}
\end{equ}
One then has the following.

\begin{lemma}\label{lem:Noise}
The quadruple $\sml(\CH_H, \{\CP_t\}_{t\ge 0}, \Wien, \{\theta_t\}_{t\ge 0}\smr)$ is a stationary noise process.
\end{lemma}
\begin{proof}
We already know from \lem{lem:separ} that $\CH_H$ is Polish. Furthermore, one has $\theta_t \circ P_t(w,\cdot\,) = w$, so
it remains to show that $\CP_t$ is a Feller transition semigroup with $\Wien$ as its unique invariant measure.
It is straightforward to check that it is a transition semigroup and the Feller property follows from the continuity of $P_t(w,\tilde w)$
with respect to $w$. By the definition \eref{e:defTrans} and the time-reversal invariance of the Wiener process, every
invariant measure for $\{\CP_t\}_{t\ge 0}$ must have its finite-dimensional distributions coincide with those of $\Wien$. Since the Borel
$\sigma$-field on $\CH_H$ is generated by the evaluation functions, this shows that $\Wien$ is the only invariant
measure.
\end{proof}

We now construct a SDS over $n$ copies of the above noise process. With a slight abuse of notation, we denote 
that noise process by $\sml(\Noise, \{\CP_t\}_{t\ge 0}, \Wien, \{\theta_t\}_{t\ge 0}\smr)$. We define the (continuous)
shift operator $R_T\from\CC \sml((-\infty,0],\R^n\smr)\to \CC_0\sml([0,T],\R^n\smr)$ by 
$\sml(R_T b\smr)(t) = b(t-T)-b(-T)$ and set
\begin{equa}[e:defPhi]
\phi \from \R_+\times \R^n \times \Noise &\to \R^n\\
	(t,x,w) &\mapsto \Phi_t\bigl(x,  R_t \CD_H w\bigr)(t)\;.
\end{equa}
From the above results, the following is straightforward:
\begin{proposition}\label{prop:defSDS}
The function $\phi$ of \eref{e:defPhi} defines a continuous SDS over the noise process 
$\sml(\Noise, \{\CP_t\}_{t\ge 0}, \Wien, \{\theta_t\}_{t\ge 0}\smr)$. Furthermore, for every generalised
initial condition $\mu$, the process generated by $\phi$ from $\mu$ is a solution to \eref{e:mainequ}
in the sense of Definition~\ref{def:solution}.
\end{proposition}
\begin{proof}
The regularity properties of $\phi$ have already been shown in \lem{lem:existSol}. The cocycle property
is an immediate consequence of the composition property for solutions of ODEs. The fact that the
processes generated by $\phi$ are solutions to \eref{e:mainequ} is a direct consequence of \eref{e:defSol}, combined
with \lem{lem:fbm}, the definition of $\CD_H$, and the fact that $\Wien$ is the Wiener measure.
\end{proof}

To conclude this section, we show that, thanks to the dissipativity condition imposed on the drift term $f$,
the SDS defined above admits any power of the Euclidean norm on $\R^n$ as a Lyapunov function:
\begin{proposition}\label{prop:Lyap}
Let $\phi$ be the continuous SDS defined above and assume that \textbf{A1} and \textbf{A2} hold. Then,
for every $p \ge 2$, the map $x \mapsto \|x\|^p$ is a Lyapunov function for $\phi$.
\end{proposition}
\begin{proof}
Fix $p\ge 2$ and let $\mu$ be an arbitrary generalised initial condition satisfying 
\begin{equ}
\int_{\R^n} \|x\|^p\,\sml(\Pi_{\R^n}^*\mu\smr)(dx) < \infty\;.
\end{equ}
Let $\tilde \phi$ be the continuous SDS 
associated by \prop{prop:defSDS} to the equation 
\begin{equ}[e:defy]
dy(t) = - y\,dt + \sigma\,dB_H(t)\;.
\end{equ}
Notice that both $\phi$ and $\tilde\phi$ are defined over the same stationary noise process.

We define $x_t$ as the process generated by $\phi$ from $\mu$ and $y_t$ as the process generated by $\tilde\phi$
from $\delta_0 \times \Wien$ (in other words $y_0 = 0$). Since both SDS are defined over the same 
stationary noise process, $x_t$ and $y_t$ are defined
over the same probability space.
The process $y_t$ is obviously Gaussian, and a direct
(but lengthy) calculation shows that its variance is given by:
\begin{equ}
\expect \|y_t\|^2 = 2H \tr \sml(\sigma \sigma^*\sml)\, e^{- t} \int_0^t s^{2H-1}\cosh(t-s)\,ds\;,
\end{equ}
In particular, one has for all $t$:
\begin{equ}[e:boundy]
\expect \|y_t\|^2 \le 2H \tr \sml(\sigma \sigma^*\sml)\, \int_0^\infty s^{2H-1}e^{-s}\,ds = \Gamma(2H+1)\tr \sml(\sigma \sigma^*\sml) \equiv C_\infty\;.
\end{equ}
Now define $z_t = x_t - y_t$. The process $z_t$ is seen to satisfy the random differential equation given by
\begin{equ}
{dz_t \over dt} = f\sml(z_t + y_t\smr) + y_t\;,\qquad z_0 = x_0\;.
\end{equ}
Furthermore, one has the following equation for $\|z_t\|^2$:
\begin{equ}
{d\|z_t\|^2 \over dt} = 2\scal{z_t, f\sml(z_t + y_t\smr)} + 2\scal{z_t,y_t}\;.
\end{equ}
Using \textbf{A2--A3} and the Cauchy-Schwarz inequality, we can estimate the right-hand side of this expression by:
\begin{equ}[e:estz]
{d\|z_t\|^2 \over dt} \le 2\C1 - 2\C2\|z_t\|^2 + 2 \scal{z_t,y_t + f(y_t)} \le -2\C2\|z_t\|^2 + \tilde C\sml(1+\|y_t\|^2\smr)^N\;,
\end{equ}
for some constant $\tilde C$. Therefore
\begin{equ}
 \|z_t\|^2 \le e^{-2\C2 t}\|x_0\|^2 + \tilde C\int_0^t e^{-2\C2(t-s)}  \sml(1+\|y_s\|^2\smr)^N\, ds\;.
\end{equ}
It follows immediately from \eref{e:boundy} and the fact that $y_s$ is Gaussian with bounded covariance \eref{e:boundy} 
that there exists a constant $C_p$ such that 
\begin{equ}
\expect \|z_t\|^p \le C_p e^{-p\C2 t} \expect \|x_0\|^p + C_p\;,
\end{equ}
for all times $t>0$. Therefore \eref{e:condLyap} holds and the proof of \prop{prop:Lyap} is complete.
\end{proof}

\section{Coupling construction}
\label{sec:coupl}

We do now have the necessary formalism to study the long-time behaviour of the SDS $\phi$ we constructed 
from \eref{e:mainequ}. The main tool that will allow us to do that is the notion of self-coupling for stochastic dynamical
systems.

\subsection{Self-coupling of SDS}
\label{sec:selfcoupl}

The main goal of this paper is to show that the asymptotic behaviour of the solutions of \eref{e:mainequ} does not
depend on its initial condition. This will then imply that the dynamics converges to a stationary state (in a suitable sense).
We therefore look for a suitable way of comparing solutions to \eref{e:mainequ}. In general, two solutions starting from
different initial points in $\R^n$ and driven with the same realisation of the noise $B_H$ have no reason of getting close
to each other as time goes by. Condition \textbf{A1} indeed only ensures that they will tend to approach each other
as long as they are sufficiently far apart. This is reasonable, since by comparing only solutions driven by the same 
realisation of the noise process, one completely forgets about the randomness of the system and the ``blurring'' this
randomness induces.

It is therefore important to compare probability measures (for example on path\-space) induced by the solutions rather
than the solution themselves. More precisely, given a SDS $\phi$ and two generalised initial conditions $\mu$ and $\nu$,
we want to compare the measures $\Evol \CQ_t\mu$ and $\Evol \CQ_t\nu$ as $t$ goes to infinity. The distance we will
work with is the total variation distance, henceforth denoted by $\|\cdot\|_\TV$. We will actually use the following useful
representation of the total variation distance. Let $\Omega$ be a measurable space and let $\prob_1$ and $\prob_2$ be two
probability measures on $\Omega$. We denote by $C(\prob_1,\prob_2)$ the set of all probability measures on $\Omega\times\Omega$
which are such that their marginals on the two components are equal to $\prob_1$ and $\prob_2$ respectively. Let furthermore
$\Delta\subset \Omega\times\Omega$ denote the diagonal, \ie the set of elements of the form $(\omega,\omega)$. We then
have
\begin{equ}[e:charTV]
\|\prob_1-\prob_2\|_\TV = 2 - \sup_{\prob \in C(\prob_1,\prob_2)} 2\prob(\Delta)\;.
\end{equ}
Elements of $C(\prob_1,\prob_2)$ will be referred to as {\em couplings} between $\prob_1$ and $\prob_2$. This leads naturally
to the following definition:

\begin{definition}\label{def:selfcoupl}
Let $\phi$ be a SDS with state space $\State$ and let $\Meas$ be the associated space of generalised initial conditions. A {\em self-coupling} 
for $\phi$ is a
measurable map $(\mu,\nu) \mapsto \Evol(\mu,\nu)$ from $\Meas\times\Meas$ into $\CD(\R_+,\State) \times \CD(\R_+,\State)$,
with the property that for every pair $(\mu,\nu)$, $\Evol(\mu,\nu)$ is a coupling for $\Evol\mu$ and $\Evol\nu$.
\end{definition}

Define the shift map $\Sigma_t \from \CD(\R_+,\State)\to\CD(\R_+,\State)$ by
\begin{equ}
\sml(\Sigma_t x\smr)(s) = x(t+s)\;.
\end{equ}
It follows immediately from the cocycle property and the stationarity of the noise process 
that $\Evol \CQ_t\mu = \Sigma_t^*\Evol\mu$. Therefore, the 
measure $\Sigma_t^*\Evol(\mu,\nu)$ is a coupling for $\Evol \CQ_t\mu$ and $\Evol \CQ_t\nu$ (which is in general
different from the coupling $\Evol(\CQ_t\mu,\CQ_t\nu)$). Our aim in the remainder of this paper is to construct a self-coupling
$\Evol(\mu,\nu)$ for the SDS associated to \eref{e:mainequ} which has the property that
\begin{equ}
\lim_{t\to\infty}\sml(\Sigma_t^*\Evol(\mu,\nu)\smr)(\Delta) = 1\;,
\end{equ}
where $\Delta$ denotes as before the diagonal of the space $\CD(\R_+,\State) \times \CD(\R_+,\State)$. 
We will then use the inequality
\begin{equ}[e:TVprop]
\|\Evol \CQ_t \mu-\Evol\CQ_t  \nu \|_\TV \le 2 - 2\sml(\Sigma_t^*\Evol(\mu,\nu)\smr)(\Delta)\;,
\end{equ}
to deduce the uniqueness of the stationary state for \eref{e:mainequ}.

In the remainder of the paper, the general way of constructing such a self-coupling will be the following. First, we fix a Polish space 
$\CA$ that contains some auxiliary information on the dynamics of the coupled process we want to keep track of. We also define a
``future'' noise space $\Noise_+$ to be equal to $\tilde \CH_H^n$, where $\tilde \CH_H$ is as in \eref{e:defTrans}. There is a natural
continuous time-shift operator on $\R \times \Noise \times \Noise_+$ defined for $t>0$ by
\begin{equ}[e:timeShift]
(s,w,\tilde w) \mapsto (s-t,P_t(w,\tilde w), S_t \tilde w)\;,\qquad \sml(S_t \tilde w\smr)(r) = \tilde w(r+t) - \tilde w(t)\;,
\end{equ}
where $P_t$ was defined in  \eref{e:defTrans}. We then construct a (measurable) map 
\begin{equa}[e:notCoupl]
\Coupl \from\State^2\times\Noise^2\times\CA &\to \R \times \pMeas(\CA \times \Noise_+^2)\;, \\
(x,y,w_x,w_y,a) &\mapsto \sml(T(x,y,w_x,w_y,a), \Wien_2(x,y,w_x,w_y,a)\smr)\;,
\end{equa}
with the properties that, for all $(x,y,w_x,w_y,a)$,
\begin{claim}
\item[\maketag{C1}{e:C1}] The time $T(x,y,w_x,w_y,a)$ is positive and greater than $1$.
\item[\maketag{C2}{e:C2}] The marginals of $\Wien_2(x,y,w_x,w_y,a)$ onto the two copies of $\Noise_+$ 
	are both equal to the Wiener measure $\Wien$.
\end{claim}
We call the map $\Coupl$ the ``coupling map'', since it yields a natural way of constructing a self-coupling for the SDS $\phi$.
The remainder of this subsection explains how to achieve this.

Given the map $\Coupl$, we can construct a Markov process on the augmented space $\BigSpace = \State^2\times \Noise^2 \times \R_+ 
\times \CA\times  \Noise_+^2 $ in the following way. As long as the component $\tau\in\R_+$ is positive, we just time-shift the
elements in $\Noise^2 \times \Noise_+^2 \times \R_+$ according to \eref{e:timeShift} and we evolve in $\CX^2$ by solving
\eref{e:mainequ}. As soon as $\tau$ becomes $0$, we redraw the future of the noise up to time $T(x,y,a)$ according to the
distribution $\Wien_2$, which may at the same time modify the information stored in $\CA$.

To shorten notations, we denote elements of $\BigSpace$ by 
\begin{equ}
X = (x,y,w_x,w_y,\tau,a,\tilde w_x, \tilde w_y)\;.
\end{equ}
With this notation, the transition function $\tilde \CQ_t$ for the process we just described is defined by:
\begin{claim}
\item[$\bullet$] For $t < \tau$, we define $\tilde \CQ_t(X;\cdot\,)$ by
\begin{equa}
\tilde \CQ_t(X;\cdot\,) = &\delta_{\phi_t(x,P_t(w_x,\tilde w_x))} \times \delta_{\phi_t(y,P_t(w_y,\tilde w_y))} \times \delta_{P_t(w_x,\tilde w_x)} \\
 &\quad   \times \delta_{P_t(w_y,\tilde w_y)} \times \delta_{\tau-t} \times \delta_{a} \times \delta_{S_t \tilde w_x}\times \delta_{S_t \tilde w_y}\;.
\end{equa}
\item[$\bullet$] For $t = \tau$, we define $\tilde \CQ_t(X;\cdot\,)$ by
\begin{equa}[e:defQtjump]
\tilde \CQ_t(X;\cdot\,) = &\delta_{\phi_t(x,P_t(w_x,\tilde w_x))} \times \delta_{\phi_t(y,P_t(w_y,\tilde w_y))} \times \delta_{P_t(w_x,\tilde w_x)}  \\
 &\quad\times \delta_{P_t(w_y,\tilde w_y)} \times \delta_{T(x,y,P_t(w_x,\tilde w_x),P_t(w_y,\tilde w_y),a)} \\
 &\quad\times \Wien_2(x,y,P_t(w_x,\tilde w_x),P_t(w_y,\tilde w_y),a) \;.
\end{equa}
\item[$\bullet$] For $t > \tau$, we define $\tilde Q_t$ by imposing that the Chapman-Kolmogorov equations hold. 
Since we assumed that $T(x,y,w_x,w_y,a)$ is always greater than $1$, this procedure is well-defined.
\end{claim}
We now construct an initial condition for this process, given two generalised initial conditions $\mu_1$ and $\mu_2$ for $\phi$.
We do this in such a way that, in the beginning, the noise component of our process lives on the diagonal of the space $\Noise^2$. In other
words, the two copies of the two-sided fBm driving our coupled system have the same past. This
is possible since the marginals of $\mu_1$ and $\mu_2$ on $\Noise$ coincide. Concerning the components of the initial condition in
$\R_+\times\CA\times\Noise_+^2$, we just draw them according to the map $\Coupl$, with some distinguished element $a_0\in\CA$.

We call $\Evol_0(\mu_1,\mu_2)$ the measure on $\BigSpace$ constructed by this procedure. Consider a cylindrical subset of $\BigSpace$
of the form
\begin{equ}
X = X_1\times X_2 \times W_1 \times W_2 \times F\;,
\end{equ}
where $F$ is a measurable subset of $\R_+ \times \CA \times \Noise_+^2$.
We make use of the disintegration $w \mapsto \mu_i^w$, yielding formally $\mu_i(dx,dw) = \mu_i^w(dx)\,\Wien(dw)$, and we
define $\Evol_0(\mu_1,\mu_2)$ by
\begin{equs}
\Evol_0(\mu_1,\mu_2)(X) = \int_{W_1\cap W_2} \int_{X_1}\int_{X_2} 
	& \sml(\delta_{T(x_1,x_2,w,w,a_0)} \times \Wien_2(x_1,x_2,w,w,a_0)\smr)(F) \\
	&\qquad \mu_2^w(dx_2)\, \mu_1^w(dx_1)\,\Wien(dw)\;.\label{e:definit}
\end{equs}
With this definition, we finally construct the self-coupling $\Evol(\mu_1,\mu_2)$ of $\phi$ corresponding to the function $\Coupl$ 
as the marginal on $\CC(\R_+,\State) \times \CC(\R_+,\State)$ of the process generated by the initial condition $\Evol_0(\mu_1,\mu_2)$
evolving under the semigroup given by $\tilde \CQ_t$. Condition \eref{e:C2} ensures that this is indeed a coupling for $\Evol\mu_1$ and 
$\Evol\mu_2$.

The following subsection gives an overview of the way the coupling function $\Coupl$ is constructed.

\subsection{Construction of the coupling function}
\label{sec:constrcoupl}

Let us consider that the initial conditions $\mu_1$ and $\mu_2$ are fixed once and for all and denote by $x_t$ and $y_t$
the two $\State$-valued processes obtained by considering the marginals of $\Evol(\mu_1,\mu_2)$
on its two $\State$ components. Define the random (but not stopping) time $\tau_\infty$ by
\begin{equ}
\tau_\infty = \inf \bigl\{t>0\,|\, \text{$x_s = y_s$ for all $s>t$}\bigr\}\;.
\end{equ}
Our aim is to find a space $\CA$ and a function $\Coupl$ satisfying \eref{e:C1} and \eref{e:C2} such that
the processes $x_t$ and $y_t$ eventually meet and stay together for all times, \ie such that
$\lim_{T\to \infty} \prob(\tau_\infty < T) = 1$. If the noise process driving the system
was Markov, the ``stay together'' part of this statement would not be a problem, since it would suffice to start driving $x_t$
and $y_t$
with identical realisations of the noise as soon as they meet. Since the fBm is not Markov, it is possible to make the future
realisations of two copies coincide with probability $1$ only if the past realisations also coincide. 
If the past realisations do not coincide for
some time, we interpret this as introducing a ``cost'' into the system, which we need to master
(this notion of cost will be made precise in Definition~\ref{def:cost} below). Fortunately, the memory
of past events becomes smaller and smaller as time goes by, which can be interpreted as a natural tendency of the
cost to decrease. This way of interpreting our system leads to the following algorithm that should be implemented by the
coupling function $\Coupl$.
\begin{equ}[e:algo]
\mhpastefig[7/8]{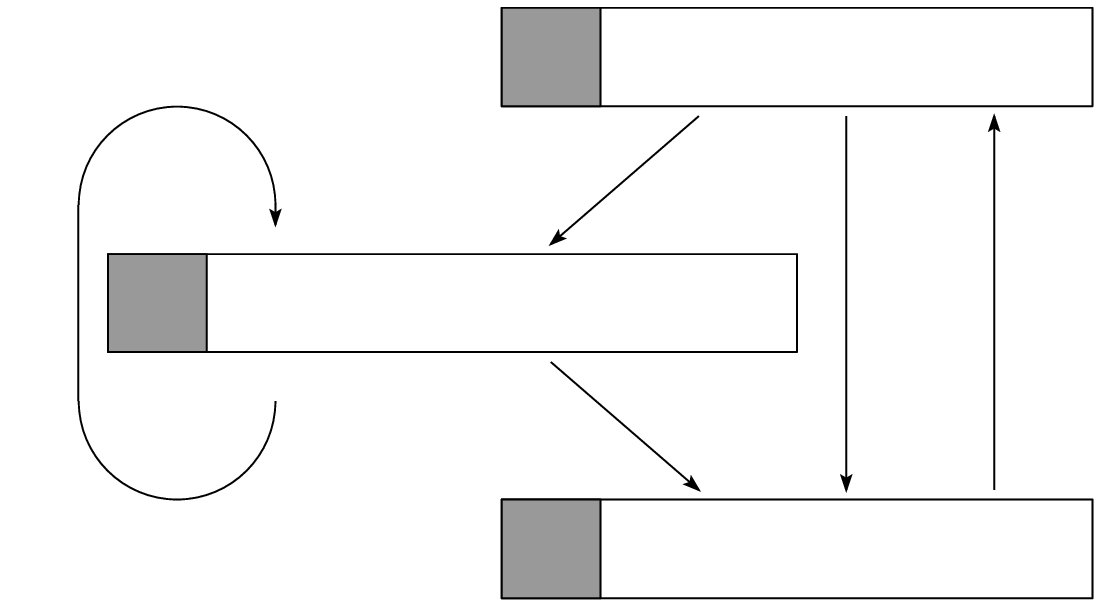}
\end{equ}
The precise meaning of the statements appearing in this diagram will be made clear in the sequel, but the general
idea of the construction should be clear by now. One step in \eref{e:algo} corresponds to the time between two
jumps of the $\tau$-component of the coupled process. Our aim is to construct the coupling function $\Coupl$ in such
a way that, with probability $1$, there is a time after which step \step{2} always succeeds. This time is then precisely the
random time $\tau_\infty$ we want to estimate.

It is clear from what has just been exposed that 
we will actually never need to consider the continuous-time process on the space $\BigSpace$ given by the
self-coupling described in the previous section, but it is sufficient to describe what happens at the beginning
of each step in \eref{e:algo}. We will therefore only consider the discrete-time dynamic obtained by sampling
the continuous-time system just before each step.
The discrete-time dynamic will
take place on the space $\CZ = \bigl(\State^2 \times \Noise^2 \times \CA\bigr) \times \R_+$ and we will
denote its elements by
\begin{equ}
(Z,\tau)\;,\quad Z = (x,y,w_x,w_y,a)\;,\quad \tau \in \R_+\;.
\end{equ}
Since the time steps of the discrete dynamic are not equally spaced, the time $\tau$ is required to keep 
track of how much time really elapsed.
The dynamic of the discrete process $(Z_n,\tau_n)$ on $\CZ$ is determined by the 
function $\Phi\from \R_+\times \CZ \times (\CA\times\Noise_+^2) \to \CZ$ given by
\begin{equs}
\Phi\bigl(t,(Z,\tau), (\tilde w_x, \tilde w_y, \tilde a)\bigr) &= \bigl( \phi_t(x,P_t(w_x,\tilde w_x)), \phi_t(y,P_t(w_y,\tilde w_y)), \\
&\quad P_t(w_x,\tilde w_x),P_t(w_y,\tilde w_y), \tilde a, \tau + t\bigr)\;.
\end{equs}
(The notations are the same as in the definition of $\tilde \CQ_t$ above.)
With this definition at hand, the transition function for the process $(Z_n,\tau_n)$ is given
by
\begin{equ}[e:transZn]
{\cscr P} (Z,\tau) = \Phi\sml(T(Z),(Z,\tau), \cdot \smr)^*\Wien_2(Z)\;,
\end{equ}
where $T$ and $\Wien_2$ are defined in \eref{e:notCoupl}. Given two generalised
initial conditions $\mu_1$ and $\mu_2$ for the original SDS, the initial condition $(Z_0,\tau_0)$ is
constructed by choosing $\tau_0 = 0$ and by drawing $Z_0$ according to the measure
\begin{equ}
\mu_0(X) = \delta_{a_0}(A)\int_{W_1\cap W_2} \int_{X_1}\int_{X_2} \mu_2^w(dx_2)\, \mu_1^w(dx_1)\,\Wien(dw)\;,
\end{equ}
where $X$ is a cylindrical set of the form $X = X_1\times X_2 \times W_1\times W_2\times A$.
It follows from the definitions \eref{e:defQtjump} and \eref{e:definit} that if we define $\tau_n$ as the $n$th jump of the process
on $\BigSpace$ constructed above and $Z_n$ as (the component in $\State^2 \times \Noise^2 \times \CA$ of) its left-hand limit at $\tau_n$, the process we obtain
is equal in law to the Markov chain that we just constructed.

Before carrying further on with the construction of $\Coupl$, we make a few preliminary computations to see how
changes in the past of the fBm affect its future. The formulae and estimates obtained in the next subsection are
crucial for the construction of $\Coupl$ and for the obtention of the bounds that lead to Theorems~\ref{theo:smallH} and \ref{theo:largeH}. In particular, \prop{prop:pastfuture} is the main estimate that leads to the coherence of the
coupling construction and to the bounds on the convergence rate towards the stationary state.

\subsection{Influence of the past on the future}

Let $w_x \in \CH_H$ and set $B_x = \CD_H w_x$. Consider furthermore two functions $g_w$ and 
$g_B$ satisfying
\begin{equ}[e:condg]
t \mapsto \int_0^t g_w(s)\,ds \in \CH_H\;,\qquad
t  \mapsto \int_0^t g_B(s)\,ds \in \CH_{1-H}\;,
\end{equ}
and define $B_y$ and $w_y$ by $B_y(0) = w_y(0) = 0$ and
\begin{equ}[e:relBw]
dB_y = dB_x + g_B\,dt\;,\qquad dw_y = dw_x + g_w\,dt\;.
\end{equ}
As an immediate consequence of the definition of $\CD_H$,
the following relations between $g_w$ and $g_B$ will ensure that $B_y = \CD_H w_y$.

\begin{lemma}\label{lem:formulas}
Let $B_x$, $B_y$, $w_x$, $w_y$, $g_B$, and $g_w$ be as in \eref{e:condg}, \eref{e:relBw} and assume that
$B_x = \CD_H w_x$ and $B_y = \CD_H w_y$. Then, $g_w$ and $g_B$ satisfy the following relation:
\minilab{e:relG}
\begin{equs}
g_w(t) &= \alpha_H {{d \over dt} \int_{-\infty}^t  (t-s)^{{1\over 2}-H}} g_B(s)\,ds\;, \label{e:relG1}\\
g_B(t) &= \gamma_H \alpha_{1-H} {{d \over dt} \int_{-\infty}^t  \Kern(t-s) g_w(s) \,ds} \;.\label{e:relG2}
\end{equs}
If $g_w(t) = 0$ for $t > t_0$, one has
\minilab{e:relG}
\begin{equ}[e:relG3]
g_B(t) = \sml(H-{\textstyle{1\over 2}}\smr)\gamma_H \alpha_{1-H}\int_{-\infty}^{t_0} (t-s)^{H - {3\over 2}} g_w(s)\,ds \;,
\end{equ}
for $t \ge t_0$. Similarly, if $g_B(t) = 0$ for $t > t_0$, one has
\minilab{e:relG}
\begin{equ}[e:relG4]
g_w(t) =\sml({\textstyle{1\over 2}}-H\smr)\alpha_{H} \int_{-\infty}^{t_0} (t-s)^{-H - {1\over 2}} g_B(s) \,ds \;,
\end{equ}
for $t \ge t_0$. If $g_w$ is differentiable for $t>t_0$ and $g_w(t) = 0$ for $t < t_0$, one has
\minilab{e:relG}
\begin{equ}[e:relG5]
g_B(t) = {\gamma_H \alpha_{1-H} g_w(t_0) \over (t-t_0)^{{1\over 2}-H}} + \gamma_H \alpha_{1-H} \int_{t_0}^t {g_w'(s) \over (t-s)^{{1\over 2}-H}}\,ds \;,
\end{equ}
for $t \ge t_0$. Similarly, if $g_B$ is differentiable for $t>t_0$ and $g_B(t) = 0$ for $t < t_0$, one has
\minilab{e:relG}
\begin{equ}[e:relG6]
g_w(t) = {\alpha_{H} g_B(t_0) \over (t-t_0)^{H-{1\over 2}}} + \alpha_{H} \int_{t_0}^t {g_B'(s) \over (t-s)^{H-{1\over 2}}}\,ds \;,
\end{equ}
for $t \ge t_0$.
\end{lemma}

\begin{proof}
The claims \eref{e:relG1} and \eref{e:relG2} follow immediately from \eref{e:relBw}, using the linearity of $\CD_H$ and the inversion formula.
The other claims are simply obtained by differentiating under the integral, see \cite{Frac} for a justification.
\end{proof}

We will be led in the sequel to consider the following situation, where $t_1$, $t_2$ and $g_1$ are assumed to be given:
\begin{equ}[e:defg]
\mhpastefig[6/7]{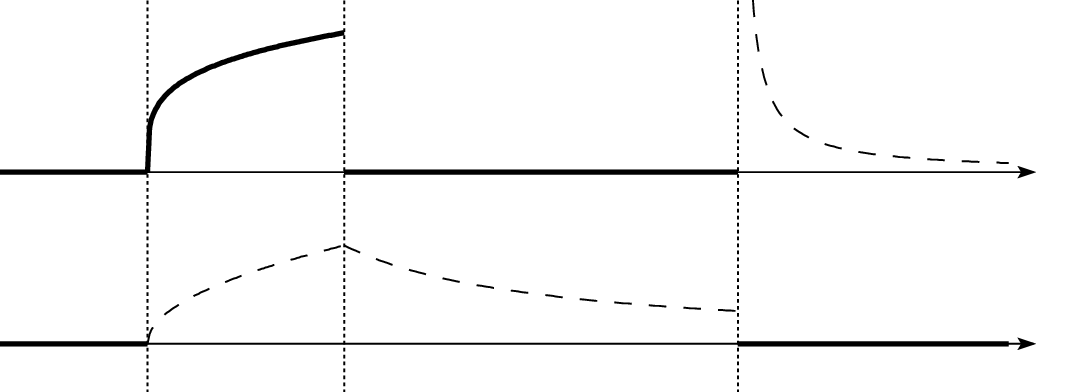}
\end{equ}
In this picture, $g_w$ and $g_B$ are related by (\ref{e:relG1}--\ref{e:relG2}) as before. The boldfaced regions indicate
that we consider the corresponding parts of $g_w$ or $g_B$ to be given. The dashed regions indicate that those
parts of $g_w$ and $g_B$ are computed from the boldfaced regions by using the relations (\ref{e:relG1}--\ref{e:relG2}).
The picture is coherent since the formulae (\ref{e:relG1}--\ref{e:relG2}) in both cases only use information about the past to compute the present.
One should think of the interval $[0,t_1]$ as representing the time spent on steps \step{1} and \step{2} of the algorithm
\eref{e:algo}. The interval $[t_1, t_2]$ corresponds to the waiting time, \ie step \step{3}. 
Let us first give an explicit formula for $g_2$ in terms of $g_1$:

\begin{lemma}
Consider the situation of \prop{prop:pastfuture}. Then, $g_2$ is given by
\begin{equ}[e:formg2]
g_2(t) = C\int_0^{t_1} { t^{{1\over 2}-H}(t_2-s)^{H-{1\over 2}} \over t + t_2 -s} g_1(s)\,ds\;,
\end{equ}
with a constant $C$ depending only on $H$.
\end{lemma}

\begin{proof}
We extend $g_1(t)$ to the whole real line by setting it equal to $0$ outside of $[0,t_1]$.
Using \lem{lem:formulas}, we see that, for some constant $C$ and for $t > t_2$,
\begin{equs}
g_2(t - t_2) &= C\int_0^{t_2} (t-s)^{-H-{1\over 2}}\,g_B(s)\,ds \\
&= C\int_0^{t_2} (t-s)^{-H-{1\over 2}}{d \over ds} \int_0^{s}(s-r)^{H-{1\over 2}}g_1(r)\,dr\, ds \\
&= C(t-t_2)^{-H-{1\over 2}} \int_0^{t_2} (t_2-r)^{H-{1\over 2}}g_1(r)\,dr \\
&\quad- C\sml(H+{\textstyle{1\over 2}}\smr) \int_0^{t_2} (t-s)^{-H-{3\over 2}}\int_0^{s}(s-r)^{H-{1\over 2}} g_1(r)\,dr\, ds\\
&\equiv C\int_0^{t_1} K(t,r) \, g_1(r)\,dr\;,
\end{equs}
where the integration stops at $t_1$ because $g_1$ is equal to $0$ for larger values of $t$.
The kernel $K$ is given by
\begin{equs}
K(t,r) &= (t-t_2)^{-H-{1\over 2}}  (t_2-r)^{H-{1\over 2}} - \sml(H+{\textstyle{1\over 2}}\smr) \int_r^{t_2} 
		(t-s)^{-H-{3\over 2}}(s-r)^{H-{1\over 2}}\, ds \\
	&= {(t-t_2)^{-H-{1\over 2}}(t_2-r)^{H+{1\over 2}} \over t-r}\Bigl({t-r \over t_2 -r} -1\Bigr) \\
	&= {(t-t_2)^{{1\over 2}-H}(t_2-r)^{H-{1\over 2}} \over t-r}\;,
\end{equs}
and the claim follows. 
\end{proof}

We give now estimates on $g_2$ in terms of $g_1$. To this end, given $\alpha > 0$,
we introduce the following norm on functions $g\from\R_+\to \R^n$:
\begin{equ}
\|g\|_\alpha^2 = \int_0^\infty (1+ t)^{2\alpha} \|g(t)\|^2\,dt\;.
\end{equ}
The following proposition is essential to the coherence of our coupling construction:

\begin{proposition}\label{prop:pastfuture}
Let $t_2 > 2 t_1 > 0$, let $g_1\from[0,t_1] \to \R^n$ be a square integrable function, and define $g_2\from\R_+\to\R_+$ by
\begin{equ}
g_2(t) = \int_0^{t_1} { t^{{1\over 2}-H}(t_2-s)^{H-{1\over 2}} \over t + t_2 -s} \|g_1(s)\|\,ds\;.
\end{equ}
Then, for every $\alpha$ satisfying
\begin{equ}
0 < \alpha < \min\{{\textstyle {1\over 2}}\,;\, H\}\;,
\end{equ}
there exists a constant $\kappa>0$ depending only on $\alpha$ and $H$ such that the estimate
\begin{equ}[e:boundfuture]
\|g_2\|_\alpha \le \kappa \Bigl|{t_2 \over t_1}\Bigr|^{\alpha -{1\over 2}} \|g_1\|_\alpha
\end{equ}
holds.
\end{proposition}

\begin{remark}
The important features of this proposition are that the constant $\kappa$ \textit{does not} depend on $t_1$ or $t_2$ and
that the exponent in \eref{e:boundfuture} is negative.
\end{remark}

\begin{proof}
We define $r = t_2/t_1$ to shorten notations.
Using \eref{e:formg2} and Cauchy-Schwarz, we then have
\begin{equs}
\|g_2(t)\| &\le C\|g_1\|_\alpha \sqrt{\int_0^{t_1} (1+ s)^{-2\alpha}  {(rt_1-s)^{2H-1} t^{1-2H} \over \sml(t + rt_1 -s\smr)^2}\,ds } \\ 
&= C\|g_1\|_\alpha {t_1}^{H-\alpha} t^{1-H-{1\over 2}} \sqrt{\int_0^{1} s^{-2\alpha}  {(r-s)^{2H-1}\over \sml(t + rt_1 - t_1s\smr)^2}\,ds }\\
&\le C \|g_1\|_\alpha {{t_1}^{H-\alpha} t^{{1\over 2}-H} r^{H-{1\over 2}} \over t + (r-1) t_1 }\;,
\end{equs}
where we made use of the assumptions that $2\alpha < 1$ and $r \ge 2$. Therefore, $\|g_2\|_\alpha$ is bounded by
\begin{equs}
\|g_2\|_\alpha &\le \kappa \|g_1\|_\alpha{t_1}^{H-\alpha} r^{H-{1\over 2}} \sqrt{\int_0^\infty {(1+ t)^{2\alpha} t^{1-2H} \over \sml(t + (r-1) t_1 \smr)^2}\,dt} \\
&\le \kappa \|g_1\|_\alpha  r^{\alpha - {1\over 2}} \sqrt{\int_0^\infty { t^{2\alpha} t^{1-2H} \over (t + 1)^2}\,dt} \;,
\end{equs}
for some constant $\kappa$, where the last inequality was obtained through the change of variables
$t \mapsto (r-1)t_1 t$ and used the fact that $r \ge 2$. The convergence of the integral is obtained under the condition $\alpha < H$ which
is verified by assumption,
so the proof of \prop{prop:pastfuture} is complete.
\end{proof}

We will construct our coupling function $\Coupl$ in such a way that there always exist functions $g_w$ and $g_B$ satisfying 
\eref{e:condg} and \eref{e:relBw}, where $w_x$ and $w_y$ denote the noise components of our coupling process, and $B_x$ and $B_y$
are obtained by applying the operator $\CD_H$ to them. 
We have now all the necessary ingredients for the construction of $\Coupl$.

\section{Definition of the coupling function}
\label{sec:defCoupl}

Our coupling construction depends on a parameter $\alpha < \min\{{1\over 2},H\}$ which we fix once and
for all. This parameter will then be tuned in Section~\ref{sec:mainproof}.

First of all, we define the auxiliary space $\CA$:
\begin{equ}[e:defA]
\CA = \{0,1,2,3\} \times \N \times \N \times \R_+\;.
\end{equ}
Elements of $\CA$ will be denoted by
\begin{equ}[e:notA]
a = (S ,N, \tilde N, T_3)\;.
\end{equ}
The component $S$ denotes which step of \eref{e:algo} is going to be performed next (the value $0$ will be used only for the initial
value $a_0$). The counter $N$ is incremented every time step \step{2} is performed and is reset to $0$
every time another step is performed. The counter $\tilde N$ is incremented every time step \step{1} or step 
\step{2} fails. If steps \step{1} or \step{2} fail, the time $T_3$ contains the duration of the upcoming 
step \step{3}. We take
\begin{equ}
a_0 = (0,1,1,0)
\end{equ}
as initial condition for our coupling construction.

Remember that the coupling function $\Coupl$ is a function from $\State^2\times \Noise^2 \times \CA$, representing the state of the system
at the end of a step, into $\R\times \pMeas(\CA\times\Noise_+^2)$, representing the duration and the realisation of the noise
for the next step. We now define $\Coupl$ for the four possible values of $S$.

\subsection{Initial stage ({\bfseries\itshape S = 0})}

Notice first that \textbf{A1} implies that
\begin{equ}[e:consA1]
{\scal{f(y)-f(x), y-x} \over \|y-x\|} \le \C4 - \C2 \|y-x\|\;,
\end{equ}
where we set $\C4 = \sqrt{\C1(\C2+\C3)}$.

In the beginning, we just wait until the two copies of our process are within distance $1 + \sml({\C4 / \C2}\smr)$ of each other.
If $x_t$ and $y_t$ satisfy \eref{e:mainequ} with the same realisation of the noise
process $B_H$, and $\rho_t = y_t - x_t$, we have by for $\|\rho_t\|$ the differential inequality
\begin{equ}
{d \|\rho_t\| \over dt} = {\scal{f(y_t)-f(x_t), \rho_t}\over \|\rho_t\|} \le \C4 - \C2\|\rho_t\|\;,
\end{equ}
and therefore by Gronwall's lemma
\begin{equ}
\|\rho_t\| \le \|y_0 - x_0\| e^{-\C2 t} + {\C4\over\C2}\bigl(1-e^{-\C2 t}\bigr)\;.
\end{equ}
It is enough to wait for a time $t = \sml(\log \|y_0 - x_0\| \smr)/ \C2$ to ensure that
$\|\rho_t\| \le 1 + \sml({\C4 / \C2}\smr)$, so we define the coupling function 
$\Coupl$ in this case by
\begin{equ}[e:defCoupl0]
T(Z,a_0) = \max\Bigl\{{\log \|y_0 - x_0\| \over \C2}\,,\, 1\Bigr\}\;,\qquad \Wien_2(Z,a_0) = \Delta^*\Wien\times \delta_{a'}\;,
\end{equ}
where the map $\Delta\from\Noise_+\to \Noise_+^2$ is defined by $\Delta(w) = (w,w)$
and the element $a'$ is given by
\begin{equ}
a' = (1,0,0,0)\;.
\end{equ}
In other terms, we wait until the two copies of the process are close to each other, and then we
proceed to step \step{1}.

\subsection{Waiting stage ({\bfseries\itshape S = 3})}

In this stage, both copies evolve with the same realisation of the underlying Wiener process.
Using notations \eref{e:notA} and \eref{e:notCoupl}, we therefore define the coupling function 
$\Coupl$ in this case by
\begin{equ}[e:defCoupl1]
T(Z,a) = T_3\;,\qquad \Wien_2(Z,a) = \Delta^*\Wien\times \delta_{a'}\;,
\end{equ}
where the map $\Delta$ is defined as above
and the element $a'$ is given by
\begin{equ}
a' = (1,N,\tilde N,0)\;.
\end{equ}
Notice that this definition is in accordance with \eref{e:algo}, \ie the counters $N$ and $\tilde N$
remain unchanged, the dynamic evolves for a time $T_3$ with two identical realisations of the
Wiener process (note that the realisations of the fBm driving the two copies of the system
are in general different, since the {\em pasts} of the Wiener processes may differ), and then
proceeds to step \step{1}.

\subsection{Hitting stage ({\bfseries\itshape S = 1})}

In this section, we construct and then analyse the map $\Coupl$ corresponding to the step \step{1}, 
which is the most important
step for our construction. We start with a few preliminary computations. Define $W_{1,1}$ as the space of
almost everywhere differentiable functions $g$, such that the quantity
\begin{equ}
 \|g\|_{1,1} = \int_0^1 \Bigl\|{dg_B(t) \over dt}\Bigr\|\,dt + \|g(0)\|\;,
\end{equ}
is finite.

\begin{lemma}\label{lem:bigH}
Let $g_B\from [0,1] \to \R^n$ be in $W_{1,1}$ and define $g_w$ by \eref{e:relG1} with 
$H\in \sml({\textstyle{1\over2}},1\smr)$. (The function $g_B$ is extended to $\R$ by setting 
it equal to $0$ outside of $[0,1]$ and $g_w$ is considered as a function
from $\R_+$ to $\R^n$.) Then, for every $\alpha \in (0,H)$, there exists a constant $C$ such that
\begin{equ}
\|g_w\|_\alpha \le C\|g_B\|_{1,1}\;.
\end{equ}
\end{lemma}

\begin{proof}
We first bound the $\L^2$ norm of $g_w$ on the interval $[0,2]$.
Using \eref{e:relG6}, we can bound $\|g_w(t)\|$ by
\begin{equ}
\|g_w(t)\| \le C \|g_B(0)\| t^{{1\over 2}-H} + C \int_0^t \|\dot g_B(s)\| (t-s)^{{1\over 2}-H}\,ds\;.
\end{equ}
Since $t^{{1\over 2}-H}$ is square integrable at the origin, it remains to bound the terms $I_1$ and $I_2$ given by
\begin{equs}
I_1 &= \int_0^2 \Bigl(\int_0^t (t-s)^{{1\over 2} -H} \|\dot g_B(s)\|\,ds \int_0^t (t-r)^{{1\over 2} -H} \|\dot g_B(r)\|\,dr \Bigr)\, dt\;,\\
I_2 &= \|g_B(0)\| \int_0^2 t^{{1\over 2} -H} \int_0^t (t-s)^{{1\over 2} -H} \|\dot g_B(s)\|\,ds\, dt\;,
\end{equs}
We only show how to bound $I_1$, as $I_2$ can be bounded in a similar fashion. Writing $r\vee s = \max\{r,s\}$ one has
\begin{equ}
I_1 = \int_0^1  \int_0^1 \int_{r \vee s}^2 (t-s)^{{1\over 2} -H} (t-r)^{{1\over 2} -H}\, dt  \|\dot g_B(s)\| \|\dot g_B(r)\|\,dr\,ds\;.
\end{equ}
Since 
\begin{equ}
\int_{r \vee s}^2 (t-s)^{{1\over 2} -H} (t-r)^{{1\over 2} -H}\,dt  \le \int_{r \vee s}^2 \sml(t-(r \vee s)\smr)^{1 - 2H}\,dt \le {2^{2-2H} \over 2-2H}\;,
\end{equ}
$I_1$ is bounded by  $C \|g_B\|_{1,1}^2$. 

It remains to bound the large-time tail of $g_w$. For $t \ge 2$, one has, again by \lem{lem:formulas},
\begin{equ}[e:largeT]
\|g_w(t)\| \le (t-1)^{-H-{1\over 2}} \sup_{s \in [0,1]} \|g_B(s)\| \le C (t-1)^{-H-{1\over 2}} \|g_B\|_{1,1}\;.
\end{equ}
It follows from the definition that the $\|\cdot\|_\alpha$-norm of this function is bounded if $\alpha < H$.
The proof of \lem{lem:bigH} is complete.
\end{proof}

In the case $H < {1\over 2}$, one has a similar result, but the regularity of $g_B$ can be weakened.

\begin{lemma}\label{lem:smallH}
Let $g_B\from [0,1] \to \R^n$ be a continuous function and define $g_w$ as in \lem{lem:bigH}, but with $H\in \sml(0,{\textstyle{1\over2}}\smr)$. Then, for every $\alpha \in (0,H)$, there exists a constant $C$ such that
\begin{equ}
\|g_w\|_\alpha \le C\sup_{t \in [0,1]} \|g_B(t)\|\;.
\end{equ}
\end{lemma}

\begin{proof}
Since $H < {1\over 2}$, one can move the derivative under the integral of the first equation in \lem{lem:formulas} to get 
\begin{equ}
\|g_w(t)\| \le C \int_0^t (t-s)^{-H -{1\over 2}}\|g_B(s)\|\,ds \le C\sup_{t \in [0,1]} \|g_B(t)\|\;.
\end{equ}
This shows that the restriction of $g_w$ to $[0,2]$ is square integrable. The large-time tail can be bounded by \eref{e:largeT}
as before.
\end{proof}

We already hinted several times towards the notion of a ``cost function'' that measures the difficulty of coupling
the two copies of the process. This notion is now made precise. Denote by $Z = (x_0,y_0,w_x,w_y)$
an element of $\State^2\times\Noise^2$ and assume that there exists a square integrable function 
$g_w\from\R_-\to\R^n$ such that
\begin{equ}[e:defg_w]
w_y(t) = w_x(t) + \int_t^0 g_w(s)\,ds\;,\qquad \forall\,t<0\;.
\end{equ}
In regard of \eref{e:formg2}, we introduce for $T>0$ the operator $\CR_T$ given by
\begin{equ}
\sml(\CR_T g\smr)(t) = C\int_{-\infty}^0 {t^{{1\over 2}-H} (T-s)^{H-{1\over 2}} \over t+T-s} \|g(s)\|\,ds\;,
\end{equ}
where $C$ is the constant appearing in \eref{e:formg2}.
The cost is then defined as follows.
\begin{definition}\label{def:cost}
The cost function $\CK_\alpha \from \L^2(\R_-) \to [0,\infty]$ is defined by
\begin{equ}[e:defcost]
\CK_\alpha (g) = \sup_{T > 0} \|\CR_T g\|_\alpha + C_K\int_{-\infty}^0 (-s)^{H-{3\over 2}} \|g(s)\|\,ds\;,
\end{equ}
where, for convenience, we define $C_K =  \sml| (2H-1)\gamma_H\alpha_{1-H}\smr|$.
Given $Z$ as above, $\CK_\alpha(Z)$ is defined as $\CK_\alpha(g_w)$ if there
exists a square integrable function $g_w$ satisfying \eref{e:defg_w} and as $\infty$ otherwise.
\end{definition}

\begin{remark}\label{rem:cost}
The cost function $\CK_\alpha$ defined above has the important property that
\begin{equ}[e:shiftcost]
\CK_\alpha(\theta_t g) \le \CK_\alpha(g)\;,\quad \text{for all $t\ge 0$,}
\end{equ}
where the shifted function $\theta_t g$ is given by
\begin{equ}
\sml(\theta_t g\smr)(s) = \cases{g(s+t) & if $s < -t$,\cr 0&otherwise. }
\end{equ}
Furthermore, it is a norm, and thus satisfies the triangle inequality.
\end{remark}

\begin{remark}\label{rem:meancost}
By \eref{e:formg2}, the first term in \eref{e:defcost} measures by how much the
two realisations of the Wiener process have to differ in order to obtain identical increments
for the associated fractional Brownian motions. By \eref{e:relG3},
the second term in \eref{e:defcost} measures by how much the two realisations
of the fBm differ if one lets the system evolve with two identical realisations of the Wiener process.
\end{remark}

We now turn to the construction of the process $(x_t,y_t)$ during step \step{1}. We will set up our coupling construction in such a way that, whenever step
\step{1} is to be performed, the initial condition $Z$ is admissible in the following sense:

\begin{definition}\label{def:defAdm}
Let $\alpha$ satisfy $0<\alpha<\min\{{1\over 2};H\}$. We say that $Z = (x_0,y_0,w_x,w_y)$ is {\em admissible} if one has
\begin{equ}[e:defAdm]
\|x_0 - y_0\| \le 1+{1 + \C4 \over \C2}\;,
\end{equ}
(the constants $\C{i}$ are as in \textbf{A1} and in \eref{e:consA1}), and its cost satisfies $\CK_\alpha(Z) \le 1$.
\end{definition}

Denote now by $\Omega$ the space
of continuous functions $\omega\from [0,1] \to \R^n$ which are the restriction to $[0,1]$ of an element of $\tilde \CH_H$. 
Our aim is construct two measures $\prob_Z^1$ and $\prob_Z^2$  on $\Omega \times \Omega$
satisfying the following conditions:
\begin{claim}
\item[\textbf{B1}] The marginals of $\prob_Z^1+\prob_Z^2$ onto the two components $\Omega$ of the product space are both equal to the Wiener
measure $\Wien$.
\item[\textbf{B2}] Let ${\cscr B}_\kappa \subset \Omega\times\Omega$ denote the set of pairs $(\tilde w_x,\tilde w_y)$ such that there exists
a function $g_w\from [0,1]\to\R^n$ satisfying
\begin{equ}
\tilde w_y(t) = \tilde w_x(t) + \int_0^t g_w(s)\,ds\;,\qquad \int_0^1 \|g_w(s)\|^2\,ds \le \kappa\;.
\end{equ}
Then, there exists a value of $\kappa$ such that, for every admissible initial condition $Z_0$, we have $\prob_Z^1({\cscr B}_\kappa) + \prob_Z^2({\cscr B}_\kappa) = 1$.
\item[\textbf{B3}] Let $(x_t,y_t)$ be the process constructed by solving \eref{e:mainequ} with respective initial conditions $x_{0}$ and
$y_{0}$, and with respective noise processes $P_t(w_x,\tilde w_x)$ and $P_t(w_y,\tilde w_y)$. Then, one has $x_1 = y_1$ for $\prob_Z^1$-almost every noise $(\tilde w_x, \tilde w_y)$. Furthermore, there exists a constant $\delta>0$ such that
$\prob_Z^1(\Omega\times\Omega) \ge \delta$ for every admissible initial condition $Z$.
\end{claim}

\begin{remark}\label{rem:}
Both measures $\prob_Z^1$ and $\prob_Z^2$ can easily be extended to measures on $\Noise_+^2$ in such a way
that \textbf{B1} holds. Since the dynamic constructed from the coupling function $\Coupl$ will not depend on this
extension, we just choose one arbitrarily and denote again by  $\prob_Z^1$ and $\prob_Z^2$ the
corresponding measures on $\Noise_+^2$.
\end{remark}

Given $\prob_Z^1$ and $\prob_Z^2$, we construct the coupling function $\Coupl$ in the following way, using notations
\eref{e:notA} and \eref{e:notCoupl}:
\begin{equ}[e:defCoupl3]
T(Z,a) = 1\;,\qquad \Wien_2(Z,a) = \prob_Z^1 \times \delta_{a_1} + \prob_Z^2 \times \delta_{a_2}\;,
\end{equ}
where the two elements $a_1$ and $a_2$ are defined as
\minilab{e:defa}
\begin{equs}
a_1 &= (2,0,\tilde N,0)\;,\label{e:defa1}\\
a_2 &= (3,0,\tilde N+1,t_* \tilde N^{4/(1-2\alpha)})\;,\label{e:defa2}
\end{equs}
for some constant $t_*$ to be determined later in this section.
Notice that this definition reflects the algorithm \eref{e:algo} and the explanation following \eref{e:notA}.
The reason behind the particular choice of the waiting time in \eref{e:defa2} will become clear in \rem{rem:step1}.

The way the construction of $ \prob_Z^1$ and $ \prob_Z^2$ works is very close to the binding construction in \cite{HExp02}. The main difference is
that the construction presented in \cite{HExp02} doesn't allow to satisfy \textbf{B2} above. We will therefore introduce
a symmetrised version of the binding construction that allows to gain a better control over $g_w$. 
If $\mu_1$ and $\mu_2$ are two positive measures with densities $D_1$ and $D_2$ with respect to
some common measure $\mu$, we define the measure $\mu_1\wedge \mu_2$ by
\begin{equ}
\sml(\mu_1\wedge \mu_2\smr)(dw) = \min\{D_1(w), D_2(w)\}\,\mu(dw)\;.
\end{equ}
The key ingredient for the construction of $ \prob_Z^1$ and $ \prob_Z^2$ is the following lemma, the proof of which will be given later in this section.
\begin{lemma}\label{lem:Psi}
Let $Z = (x_0,y_0,w_x,w_y)$ be an admissible initial condition 
and let $H$, $\sigma$, and $f$ satisfy the hypotheses of either
\theo{theo:smallH} or \theo{theo:largeH}. Then, there exists a measurable map $\Psi_Z\from\Omega\to\Omega$ with
measurable inverse,
having the following properties.
\begin{claim}
\item[\textbf{B1'}] There exists a constant $\delta > 0$ such that $\Wien \wedge \Psi_Z^*\Wien$ has mass bigger than $2\delta$
for every admissible initial condition $Z$.
\item[\textbf{B2'}] There exists a constant $\kappa$ such that $\{(\tilde w_x,\tilde w_y)\,|\, \tilde w_y = \Psi_Z(\tilde w_x)\} \subset {\cscr B}_\kappa$
for every admissible initial condition $Z$.
\item[\textbf{B3'}] Let $(x_t,y_t)$ be the process constructed by solving \eref{e:mainequ} with respective initial conditions $x_{0}$ and $y_{0}$, and with noise processes $P_t(w_x,\tilde w_x)$ and $P_t(w_y,\Psi_Z(\tilde w_x))$. Then, one
has $x_1 = y_1$ for every $\tilde w_x \in \Omega$ and every admissible initial condition $Z$.
\end{claim}
Furthermore, the maps $\Psi_Z$ and $\Psi_Z^{-1}$ are measurable with respect to $Z$.
\end{lemma}

Given such a $\Psi_Z$, we first define the maps $\Psi_\uparrow$ and $\Psi_\rightarrow$
from $\Omega$ to $\Omega\times\Omega$ by
\begin{equ}
\Psi_\uparrow(\tilde w_x) = \sml(\tilde w_x,\Psi_Z(\tilde w_x)\smr)\;,\qquad\Psi_\rightarrow(\tilde w_y) = \sml(\Psi_Z^{-1}(\tilde w_y),\tilde w_y\smr)\;.
\end{equ}
(See also Figure~\ref{fig:constr} below.)
We also define the ``switch map'' $S\from\Omega\times\Omega \to \Omega\times\Omega$ by $S(\tilde w_x,\tilde w_y) = (\tilde w_y,\tilde w_x)$.
\begin{figure}[h]
\begin{center}
\mhpastefig{Coupling}
\end{center}
\caption{Construction of $\prob_Z$.}\label{fig:constr}
\end{figure}

With these definitions at hand, we construct two measures $\prob_Z^1$ and  $\tilde \prob_Z^1$ on $\Omega\times\Omega$ by
\begin{equ}[e:constP]
\prob_Z^1 = {1\over 2} \bigl(\Psi_\uparrow^*\Wien \wedge \Psi_\rightarrow^*\Wien\bigr)\;,\qquad \tilde\prob_Z^1 = \prob_Z^1 + S^*\prob_Z^1 \;.
\end{equ}
On Figure~\ref{fig:constr}, $\prob_Z^1$ lives on the boldfaced curve and $\tilde \prob_Z^1$ is its symmetrised version
which lives on both the boldfaced and the dashed curve.
Denote by $\Pi_{i}\from\Omega\times\Omega \to \Omega$ the projectors onto the $i$th component and by
$\Delta\from\Omega\to\Omega \times \Omega$ the lift onto the diagonal $\Delta(w) = (w,w)$. Then, we define
the measure $\prob_Z^2$ by
\begin{equ}[e:defP]
\prob_Z^2 = S^*\prob_Z^1 + \Delta^*\sml(\Wien - \Pi_1^* \tilde \prob_Z^1\smr)\;.
\end{equ}
By \eref{e:constP}, $\Wien > \Pi_1  \tilde \prob_Z^1$, so $\prob_Z^1$ and $\prob_Z^2$ are both positive
and their sum is a probability measure. 
Furthermore, one has by definition
\begin{equ}
\prob_Z^1 + \prob_Z^2 = \tilde \prob_Z^1+ \Delta^*\sml(\Wien - \Pi_1^* \tilde \prob_Z^1\smr)\;.
\end{equ}
Since $\Pi_1^*\Delta^*$ is the identity, this immediately implies
\begin{equ}
\Pi_1^*\prob_Z^1 + \Pi_1^*\prob_Z^2= \Wien\;.
\end{equ}
The symmetry $S^*\tilde \prob_Z^1 = \tilde \prob_Z^1$ then implies that the second marginal is also equal
to $\Wien$,
\ie \textbf{B1} is satisfied. Furthermore, the set $\{(\tilde w_x,\tilde w_y)\,|\, \tilde w_y = \Psi_Z(\tilde w_x)\}$ has $\prob_Z$-measure
bigger than $\delta$ by \textbf{B1'}, so \textbf{B3} is satisfied as well. Finally, \textbf{B2} is an immediate consequence of
\textbf{B2'}. It remains to construct the function $\Psi_Z$.

\begin{proof}[of \lem{lem:Psi}]
As previously, we write $Z$ as
\begin{equ}[e:notZ0]
Z = (x_0,y_0,w_x,w_y)\;.
\end{equ}
In order to construct $\Psi_Z$, we proceed as in \cite[Sect.\ 5]{HExp02}, except that we want the solutions $x_t$ and $y_t$
to become equal after time $1$. Let $\tilde w_x\in\Omega$ be given and define
\begin{equ}[e:defBH]
B_H(t) = \sml(\CD_H P_1(w_x,\tilde w_x)\smr)(t-1)\;,
\end{equ}
where $W$ denotes the corresponding part of the initial condition $Z_0$ in \eref{e:notZ0}. We write the solutions
to \eref{e:mainequ} as
\minilab{e:defxy}
\begin{equs}
dx_t &= f(x_t)\,dt + \sigma dB_H(t) \label{e:defxt}\;,\\
dy_t &= f(y_t)\,dt + \sigma dB_H(t) + \sigma \tilde g_B(t)\,dt\;, \label{e:defyt}
\end{equs}
where $\tilde g_B(t)$ is a function to be determined. Notice that $x_t$ is completely determined by $\tilde w_x$ and by the
initial condition $Z$. We introduce the process $\rho_t = y_t - x_t$, so we get
\begin{equ}[e:equrho]
{d \rho_t \over dt} = f(x_t+ \rho_t) - f(x_t) + \sigma \tilde g_B(t)\;.
\end{equ}
We now define $\tilde g_B(t)$ by
\begin{equ}[e:defgb]
\tilde g_B(t) = -\sigma^{-1}\Bigl(\kappa_1\rho_t + \kappa_2 {\rho_t \over \sqrt{\|\rho_t\|}}\Bigr)\;,
\end{equ}
for two constants $\kappa_1$ and $\kappa_2$ to be specified. This yields for the
norm of $\rho_t$ the estimate
\begin{equ}
{d \|\rho_t\|^2 \over dt} \le 2(\C3 - \kappa_1)\|\rho_t\|^2 - 2\kappa_2 \|\rho_t\|^{3/2}\;.
\end{equ}
We choose $\kappa_1 = \C3$ and so
\begin{equ}[e:estrho]
\|\rho_t\| \le \cases{\bigl(6\kappa_2 t - \sqrt{\|\rho_0\|}\bigr)^2 & for 
	$t < \sqrt{\|\rho_0\|} / (6\kappa_2)$,\cr 0& for $t \ge \sqrt{\|\rho_0\|} / (6\kappa_2)$.}
\end{equ}
We can then choose
$\kappa_2$ sufficiently large, so that $\|\rho_t\| = 0$ for $t > 1/2$. Since the initial condition
was admissible by assumption, the constant $\kappa_2$ can be chosen as a function of the constants $\C{i}$ only.
Notice also that the preceding construction yields $\tilde g_B$ as a function of
$Z$ and $\tilde w_x$ only.

We then construct $\tilde w_y = \Psi_Z(\tilde w_x)$ in such a way that \eref{e:defxy} is satisfied with the function
$\tilde g_B$ we just constructed. Define $g_w$ by \eref{e:defg_w} and construct $g_B$ by applying \eref{e:relG2}. Then,
we extend $\tilde g_B$ to $(-\infty,1]$ by simply putting it equal to $g_B$ on $(-\infty,0]$. Applying the inverse formula
\eref{e:relG1}, we obtain a function $\tilde g_w$ on  $(-\infty,1]$, which is equal to $g_w$ on $(-\infty,0]$ and which is
such that 
\begin{equ}
\sml(\Psi_Z(\tilde w_x)\smr)(t) \equiv \tilde w_x(t) + \int_0^t \tilde g_w(s)\,ds \;,
\end{equ}
has precisely the required property.

It remains to check that the family of maps $\Psi_Z$ constructed this way has the properties stated in \lem{lem:Psi}.
The inverse of $\Psi_Z$ is constructed in the following way. Choose $\tilde w_y \in \Omega$ and consider the solution
to the equation
\begin{equ}
dy_t = f(y_t)\,dt + \sigma dB_H'(t)\;,
\end{equ}
where $B_H$ is defined as in \eref{e:defBH} with $x$ replaced by $y$. Once $y_t$ is obtained, one can construct the
process $\rho_t$ as before, but this time by solving
\begin{equ}
{d \rho_t \over dt} = f(y_t) - f(y_t - \rho_t)  - \Bigl(\kappa_1\rho_t + \kappa_2 {\rho_t \over \sqrt{\|\rho_t\|}}\Bigr)\;.
\end{equ}
This allows to define $\tilde g_B$ as in \eref{e:defgb}. The element $\tilde w_x \equiv \Psi_Z^{-1}(\tilde w_y)$ is then obtained by the same
procedure as before.

Before turning to the proof of properties \textbf{B1'}--\textbf{B3'}, we give some estimate on the function $\tilde g_w$ that we just constructed.

\begin{lemma}\label{lem:estg_w}
Assume that the conditions of \lem{lem:Psi} hold. Then, there exists a constant $K$ such that the function $\tilde g_w(Z,\tilde w_x)$ 
constructed above satisfies
\begin{equ}
\int_0^1 \sml\|\tilde g_w(Z,\tilde w_x)(s)\smr\|^2\,ds < K\;,
\end{equ}
for every admissible initial condition $Z$ and for every $\tilde w_x \in \Noise_+$.
\end{lemma}

\begin{proof}
We write $\tilde g_w(t)$ for $t>0$ as
\begin{equs}
\tilde g_w(t) &= C\int_{-\infty}^0 {t^{{1\over 2}-H}(-s)^{H-{1\over 2}} \over t-s} g_w(s)\,ds + \alpha_H {d\over dt}\int_0^t
(t-s)^{{1\over 2}-H} \tilde g_B(s)\,ds\;, \\
&\equiv \tilde g_w^{(1)}(t) + \tilde g_w^{(2)}(t)\;.
\end{equs}
where $g_w$ is defined by \eref{e:defg_w}, $g_B$ is given by \eref{e:defgb}, and the constant $C$ is the constant appearing
in \eref{e:formg2}. 
The $\L^2$-norm of $\tilde g_w^{(1)}$ is bounded by $1$ by 
the assumption that $Z$ is admissible.
To bound the norm of $\tilde g_w^{(2)}$, we treat the cases $H<{1\over 2}$ and $H>{1\over 2}$ separately.

\noindent \textbf{The case $\mathbf{H < {1\over 2}}$.} For this case, we simply combine \lem{lem:smallH} with the definition
\eref{e:defgb} and the estimate \eref{e:estrho}. 

\noindent \textbf{The case $\mathbf{H>{1\over 2}}$.} For this case, we apply \lem{lem:bigH}, so we bound the $\|\cdot\|_{1,1}$-norm of $\tilde g_B$. By \eref{e:defgb}, one has
\begin{equ}[e:estdgb]
\Bigl\|{d \over dt} \tilde g_B(t)\Bigr\| \le C\Bigl\|{d\rho_t \over dt}\Bigr\| \bigl(1 + \|\rho_t\|^{-{1\over 2}}\bigr)\;,
\end{equ}
for some positive constant $C$. Using \eref{e:equrho}, the assumption about the boundedness
of the derivative of $f$, and the definition \eref{e:defgb} we get
\begin{equ}
 \Bigl\|{d\rho_t \over dt}\Bigr\| \le C \bigl(\|\rho_t\| + \sqrt{\|\rho_t\|}\bigr) \;.
\end{equ}
Combining this with \eref{e:estdgb} and \eref{e:estrho}, the required bound on $\|\tilde g_B\|_{1,1}$ follows.
\end{proof}

Property \textbf{B1'} now follows from \lem{lem:estg_w} and Girsanov's theorem in the following way. Denote by $\Dens_Z$ the density of $\Psi_Z^*\Wien$ with
respect to $\Wien$, \ie $\sml(\Psi_Z^*\Wien\smr)(d\tilde w_x) = \Dens_Z(\tilde w_x)\,\Wien(d\tilde w_x)$. It is given by Girsanov's formula
\begin{equ}
\Dens_Z(\tilde w_x) = \exp \Bigl(\int_0^1 \bigl\langle\sml(\tilde g_w(Z,\tilde w_x)\smr)(t)\,,\,d\tilde w_x(t)\bigr\rangle -{1\over 2}\int_0^1 \sml\|\tilde g_w(Z,\tilde w_x)\smr\|^2(t)\,dt\Bigr)\;.
\end{equ}
One can check (see \eg \cite{MatNS}) that $\|\Wien \wedge \Psi_Z^*\Wien\|_\TV$ is bounded from below by
\begin{equ}
\|\Wien \wedge \Psi_Z^*\Wien\|_\TV \ge \Bigl(4\int_\Omega \Dens_Z(w)^{-2}\,\Wien(dw) \Bigr)\;.
\end{equ}
Property \textbf{B1'} thus follows immediately from \lem{lem:estg_w}, using the fact that
\begin{equ}
\int_\Omega \exp \Bigl(-2\int_0^1 \bigl\langle\sml(\tilde g_w(Z,\tilde w_x)\smr)(t)\,,\,d\tilde w_x(t)\bigr\rangle - 2\int_0^1 \sml\|\tilde g_w(Z,\tilde w_x)\smr\|^2(t)\,dt\Bigr)\,\Wien(dw) = 1\;.
\end{equ}
 Property \textbf{B2'} is also an immediate consequence of \lem{lem:estg_w}, and property \textbf{B3'} follows by construction from \eref{e:estrho}. The proof of \lem{lem:Psi} is complete.
\end{proof}

Before concluding this subsection we show that, if step \step{1} fails, $t_*$ can be chosen in such a way that
the waiting time $t_*\tilde N^{4/(1-2\alpha)}$ in
\eref{e:defa2} is long enough so that \eref{e:defAdm} holds again after step \step{3} and so that the cost
 function does not  increase by more than $1/(2\tilde N^2)$. By the triangle inequality, the second claim follows
 if we show that
\begin{equ}[e:boundCost]
\CK_\alpha\sml(\theta_t \tilde g_w(Z,\tilde w_x)\smr) \le {1 \over 2 \tilde N^2}\;,
\end{equ}
whenever $t$ is large enough (the shift $\theta_t$ is as in \eref{e:shiftcost}). Combining \eref{e:boundfuture}, \lem{lem:estg_w}, and the
definition of $\CK_\alpha$, we get, for some constant $C$,
\begin{equ}
\CK_\alpha\sml(\theta_t \tilde g_w(Z,\tilde w_x)\smr) \le C t^{\alpha - {1\over 2}} + C t^{H-{3\over 2}}\;,\quad\text{for $t\ge 2$.}
\end{equ}
There thus exists a constant $t_*$ such that the bound \eref{e:boundCost} is satisfied if the waiting time 
is longer than $t_* \tilde N^{4/(1-2\alpha)}$. It remains to show that \eref{e:defAdm} holds after the waiting
time is over. If step \step{1} failed, the realisations $\tilde w_x$ and $\tilde w_y$ are drawn either in the set
\begin{equ}
\tilde \Delta_1 = \{ (\tilde w_x, \tilde w_y) \in \Omega^2 \,|\, \tilde w_x = \tilde w_y \}\;,
\end{equ}
or in the set
\begin{equ}
\tilde \Delta_2 = \{(\tilde w_x, \tilde w_y) \in \Omega^2 \,|\, \tilde w_x = \Psi_Z(\tilde w_y) \}
\end{equ}
(see Figure~\ref{fig:constr}). In order to describe the dynamics also during the waiting time
(\ie step \step{3}), we extend those sets to $\Noise_+^2$ by
\begin{equs}
\Delta_i &= \big\{(\tilde w_x, \tilde w_y) \in \Noise_+^2\,|\, (\tilde w_x |_{[0,1]}, \tilde w_y |_{[0,1]}) \in \tilde \Delta_i\;,
\\
&\qquad \text{and}
\quad \tilde w_x (t) - \tilde w_y (t) = \text{const}\quad \text{for $t>1$}\big\}\;.
\end{equs}
Given an admissible initial condition $Z = (x_0,y_0,w_x,w_y)$ and a pair $(\tilde w_x, \tilde w_y) \in \Noise_+^2$, we consider
the solutions $x_t$ and $y_t$ to \eref{e:mainequ} given by
\begin{equa}[e:evolCoupl]
dx_t &= f(x_t)\,dt + \sigma dB_H^x(t) \;,\\
dy_t &= f(y_t)\,dt + \sigma dB_H^y(t) \;,
\end{equa}
where $B_H^x$ (and similarly for $B_H^y$) is constructed as usual by concatenating $w_x$ and $\tilde w_x$ and applying the operator
$\CD_H$. 
The key observation is the following lemma.

\begin{lemma}\label{lem:controldist}
Let $Z$ be an admissible initial condition as above, let $(\tilde w_x, \tilde w_y) \in \Delta_1 \cup \Delta_2$, and 
let $x_t$ and $y_t$ be given by \eref{e:evolCoupl} for $t>0$. Then, there exists a constant $t_* > 0$ such that
\begin{equ}
\|x_t - y_t\| \le 1+ {1 + \C4 \over \C2}
\end{equ}
holds again for $t > t_*$.
\end{lemma}

\begin{proof}
Fix an admissible initial condition $Z$ and consider the case when $(\tilde w_x, \tilde w_y) \in \Delta_2$ first. 
Let $g_w \from \R_- \to \R^n$ be as in \eref{e:defg_w} and define $\tilde g_w \from \R_+ \to \R^n$
by
\begin{equ}
\tilde w_y(t) = \tilde w_x(t) + \int_0^t \tilde g_w(s)\,ds\;.
\end{equ}
Introducing $\rho_t = y_t - x_t$, we see that it
satisfies the equation
\begin{equ}[e:equrho2]
{d \rho_t \over dt} = f(y_t) - f(x_t) + \sigma \CG_t\;,
\end{equ}
where the function $\CG_t$ is given by
\begin{equ}[e:defG]
\CG_t = c_1 \int_{-\infty}^0 (t-s)^{H-{3\over 2}} g_w(s)\,ds + 
c_2 {d \over dt} \int_{0}^t (t-s)^{H-{1\over 2}} \tilde g_w(s)\,ds\;,
\end{equ}
with some constants $c_1$ and $c_2$ depending only on $H$.
It follows from \eref{e:equrho2}, \eref{e:consA1}, and Gronwall's lemma, that the Euclidean norm $\|\rho_t\|$ satisfies the inequality
\begin{equ}[e:estrhoGron]
\|\rho_t\| \le e^{-\C2 t}\|\rho_0\| + \int_0^t e^{-\C2 (t-s)}\sml(\C4 + \|\CG_s\|\smr)\,ds\;.
\end{equ}
Consider first the time interval $[0,1]$ and define
\begin{equ}
\tilde \CG_t = c_1\int_{-\infty}^0 (t-s)^{H-{3\over 2}} g_w(s)\,ds -
c_2 {d \over dt} \int_{0}^t (t-s)^{H-{1\over 2}} \tilde g_w(s)\,ds\;,
\end{equ}
\ie, we simply reversed the sign of $\tilde g_w$.
This corresponds to the case where $(\tilde w_x, \tilde w_y)$ are interchanged, and thus satisfy
$\tilde w_y = \Psi_Z(\tilde w_x)$ instead of $\tilde w_x = \Psi_Z(\tilde w_y)$.
We thus deduce from \eref{e:defgb} and \eref{e:estrho} that 
\begin{equ}[e:wantboundG]
\|\tilde \CG_s\| \le \|\sigma^{-1}\|\bigl(\kappa_1 \|\rho_0\|Ê+ \kappa_2 \sqrt{\|\rho_0\|}\bigr)\;,
\end{equ}
for $s \in [0,1]$. This yields for $\|\CG_s\|$ the estimate
\begin{equs}
\|\CG_s\| &\le \|\sigma^{-1}\|\bigl(\kappa_1 \|\rho_0\|Ê+ \kappa_2 \sqrt{\|\rho_0\|}\bigr) + 2 c_1 \int_{-\infty}^0 (t-s)^{H-{3\over 2}} \|g_w(s)\|\,ds \\
&\le  \|\sigma^{-1}\|\bigl(\kappa_1 \|\rho_0\|Ê+ \kappa_2 \sqrt{\|\rho_0\|}\bigr) + 1 \;, \label{e:estG1}
\end{equs}  
where we used the fact that $Z$ is admissible for the second step. Notice that \eref{e:estG1} only holds for
$s \in [0,1]$, so we consider now the case $s > 1$. In this case, we can write $\CG_t$ as
\begin{equ}
\CG_t = c_1 \int_{-\infty}^0 (t-s)^{H-{3\over 2}} g_w(s)\,ds + 
c_1 \int_{0}^1 (t-s)^{H-{3\over 2}} \tilde g_w(s)\,ds\;.
\end{equ}
The first term is bounded by $1$ as before.
In order to bound the second term, we use \lem{lem:estg_w}, so we get
\begin{equ}[e:boundGt]
\|\CG_t\| \le  1 + \sqrt{{K\over 2H-2} \bigl((t-1)^{2H-2} - t^{2H-2}\bigr)}\;.
\end{equ}
This function has a singularity at $t=1$, but this singularity is always integrable. For $t>2$ say, it behaves like
$t^{H-{3\over 2}}$. Putting the estimates \eref{e:estG1} and \eref{e:boundGt} into \eref{e:estrhoGron}, we
see that there exists a constant $C$ depending only on $H$ and on the parameters in assumption
\textbf{A1} such that, for $t>2$, one has the estimate
\begin{equ}
\|\rho_t\| \le e^{-\C2 t}\|\rho_0\| + {1+\C4 \over \C2} + C t^{H-{3\over 2}}\;.
\end{equ}
The claim follows at once.
\end{proof}

\begin{remark}\label{rem:step1}
To summarise, we have shown the following in this section:
\begin{claim}[AA]
\item[1.] There exists a positive constant $\delta$ such that if the state $Z$ of the coupled system is admissible, 
step \step{1} has a probability larger than $\delta$ to succeed. 
\item[2.] If step \step{1} fails and the waiting time for step \step{3} is chosen larger than 
$t_*\tilde N^{4/(1-2\alpha)}$, then the state of the coupled 
system is again admissible after the end of step \step{3}, provided the cost $\CK_\alpha(Z)$ at the beginning of step \step{1}
was smaller than $1-{1 \over 2\tilde N^2}$.
\item[3.] The increase in the cost given between the beginning of step \step{1} and the end
of step \step{3} is smaller than ${1 \over 2\tilde N^2}$.
\end{claim}
\end{remark}

In the following subsection, we will define step \step{2} and so conclude the construction and the analysis of the
coupling function $\Coupl$.

\subsection{Coupling stage ({\bfseries\itshape S = 2})}

In this subsection, we construct and analyse the coupling map $\Coupl$ corresponding to step \step{2}. 
Following \eref{e:algo}, we construct it in such a way that, with positive probability, the two copies of the process
\eref{e:mainequ} are driven with the same noise. In other terms, if $Z=(x_0,y_0,w_x,w_y)$ denotes the
state of our coupled system at the beginning of step \step{2}, we construct a measure $\prob_Z$
on $\Noise_+^2$ such that if $(\tilde w_x, \tilde w_y)$ is drawn according to $\prob_Z$, then one has
\begin{equ}[e:defcoupl]
\sml(\CD_H(w_x \sqcup \tilde w_x)\smr)(t) = \sml(\CD_H(w_y \sqcup \tilde w_y)\smr)(t)\;,\quad t>0\;,
\end{equ}
with positive probability. Here, $\sqcup$ denotes the concatenation operator given by
\begin{equ}
\sml(w \sqcup \tilde w\smr)(t) = \cases{w(t) & for $t < 0$,\cr\tilde w(t) & for $t \ge 0$.}
\end{equ}
In the notation \eref{e:notA}, step \step{2} will have a duration $2^N$ and $N$ will be incremented
by $1$ every time step \step{2} succeeds.

The construction of $\prob_Z$ will be similar in spirit to the construction of the previous section.
We therefore introduce as before the function $\tilde g_w$ given by
\begin{equ}[e:defgwtilde]
\tilde w_y(t) = \tilde w_x(t) + \int_0^t \tilde g_w(s)\,ds\;.
\end{equ}
Our main concern is of course to get good bounds on this function $\tilde g_w$. This is
achieved by the following lemma, which is crucial in the process of showing that step \step{2} will
eventually succeed infinitely often.

\begin{lemma}\label{lem:step1}
Let $Z_0$ be an admissible initial condition and denote by $\CT$ the measure on $\State^2\times\Noise^2$
obtained by evolving $Z_0$ according to the successful realisation of step \step{1}. Then, there exists a constant
$\tilde K > 0$ depending only on $H$, $\alpha$, and the parameters appearing in \textbf{A1}, such that for $\CT$-almost
every $Z = (x,y,w_x,w_y)$, and for every pair $(\tilde w_x,\tilde w_y)$ satisfying \eref{e:defcoupl}, we have the
bounds
\begin{equ}[e:boundstep2]
\|\tilde g_w\|_\alpha \le \tilde K \;,\qquad \Bigl\|{d \tilde g_w \over dt}\Bigr\|_{\alpha+1} \le \tilde K \;.
\end{equ}
Furthermore, one has $x=y$, $\CT$-almost surely.
\end{lemma}

\begin{proof}
It is clear from \lem{lem:Psi} that $x=y$.
Let now $Z$ be an element drawn according to $\CT$ 
and denote by $g_w\from\R_-\to \R^n$ the function formally defined by
\begin{equ}[e:defgw2]
dw_y(t)  = dw_x(t) + g_w(t)\,dt\;.
\end{equ}
We also denote by $g_b\from \R_- \to \R^n$ the function such that
\begin{equ}[e:defgb2]
dB_y(t)  = dB_x(t) + g_b(t)\,dt\;,
\end{equ} 
where $B_x = \CD_H w_x$ and $B_y = \CD_H w_y$. 
(Note that $g_w$ and $g_b$ are almost surely well-defined, so we discard elements $Z$
for which they can not be defined.)
Since $Z$ corresponds almost surely to a successful realization of step \textbf{1}, $g_b$ is equal
on the interval $[-1,0]$ (up to translation in time) to the function $\tilde g_B$ constructed in \eref{e:defgb}. By \eref{e:estrho}, there exists therefore a constant $C_g$ such that
\begin{equ}[e:boundgb]
\|g_b(s)\| \le \cases{C_g & for $s\in [-1,-{\textstyle{1\over 2}})$,\cr 0& for $s\in[-{\textstyle{1\over 2}},0]$.}
\end{equ}
Combining the linearity of $\CD_H$ with
\eref{e:formg2}, one can see that if $(\tilde w_x,\tilde w_y)$ satisfy \eref{e:defcoupl}, then the
function $\tilde g_w$ is given by the formula
\begin{equ}[e:exprgwtilde1]
\tilde g_w(t) =  C_1\int_{-\infty}^{-1} {|t+1|^{{1\over 2}-H}|s + 1|^{H-{1\over 2}} \over t-s}g_w(s)\,ds
+ C_2\int_{-1}^{-1/2} (t-s)^{-H-{1\over 2}} g_b(s)\,ds\;,
\end{equ}
for some constants $C_1$ and $C_2$ depending only on $H$. Notice that the second integral 
only goes up to $1/2$ because of \eref{e:boundgb}.

Since the initial condition $Z_0$  is admissible by assumption, the $\|\cdot\|_\alpha$-norm of the first
term is bounded by $1$. The $\|\cdot\|_\alpha$-norm of the second term is also bounded by a constant, 
using \eref{e:boundgb} and the assumption
$\alpha < H$.

Deriving \eref{e:exprgwtilde1} with respect to $t$, we see that there exists a constant $K$ such that
\begin{equa}[e:boundderg]
\Bigl\|{d\tilde g_w(t) \over dt}\Bigr\| &\le {K \over t+1} \Big(\int_{-\infty}^{-1} {|t+1|^{{1\over 2}-H}|s + 1|^{H-{1\over 2}} \over t-s} \|g_w(s)\|\,ds
\\
&\qquad\qquad+ \int_{-1}^{-1/2} (t-s)^{-H-{1\over 2}} \|g_b(s)\|\,ds\Big)\;,
\end{equa}
and the bound on the derivative follows as previously.
\end{proof}

The definition of our coupling function will be based on the following lemma:
\begin{lemma}\label{lem:gauss}
Let $\CN$ be the normal distribution on $\R$, choose $a \in \R$, $b\ge |a|$, and define 
$M = \max\{4b,2\log(8/b)\}$. Then, there exists a measure
$\CN_{a,b}^2$ on $\R^2$ satisfying the following properties:
\begin{claim}[33]
\item[1.] Both marginals of $\CN_{a,b}^2$ are equal to $\CN$.
\item[2.] If $|b|\le 1$, one has
\begin{equ}
\CN_{a,b}^2 \bigl(\bigl\{(x,y)\,|\, y = x+a\bigr\}\bigr) > 1 - b\;.
\end{equ}
Furthermore, the above quantity is always positive.
\item[3.] One has
\begin{equ}
\CN_{a,b}^2 \bigl(\bigl\{(x,y)\,|\, |y - x| \le M\bigr\}\bigr) = 1\;.
\end{equ}
\end{claim}
\end{lemma}

\begin{proof}
Consider the following picture:
\begin{center}
\begin{minipage}{5cm}
\mhpastefig{Gaussian}
\end{minipage}
\begin{minipage}{4cm}
\begin{equs}
L_1\,:\quad y &= x \;,\\
L_2\,:\quad y &= -x \;,\\
L_3\,:\quad y &= x+a \;.\\
\quad&\quad
\end{equs}
\end{minipage}
\end{center}
Denote by $\CN_x$ the normal distribution on the set $L_x = \{(x,y)\,|\, y = 0\}$ and by $\CN_y$ the normal
distribution on the set  $L_y = \{(x,y)\,|\, x = 0\}$. We also define the maps $\pi_{i,x}$ (respect. $\pi_{i,y}$) from
$L_x$ (respect. $L_y$) to $L_i$, obtained by only modifying the $y$ (respect. $x$) coordinate. Notice that 
these maps are invertible and denote their inverses by $\tilde \pi_{i,x}$ (respect. $\tilde \pi_{i,y}$).
We also denote by $\CN_x|_M$ (respect. $\CN_y|_M$) the restriction of $\CN_x$ (respect. $\CN_y$)
to the square $[-{M\over 2},{M\over 2}]^2$.

With these notations, we define the measure $\CN_3$ on $L_3$ as
\begin{equ}
\CN_3 = \pi_{3,x}^* \sml(\CN_x |_M\smr)\wedge \pi_{3,y}^* \sml(\CN_y |_M\smr)\;.
\end{equ}
The measure $\CN_{a,b}^2$ is then defined as
\begin{equ}
\CN_{a,b}^2 = \CN_3 +  \pi_{2,x}^*\bigl(\sml(\CN_x |_M\smr) - \tilde \pi_{3,x}^*\CN_3\bigr)
+ \pi_{1,x}^*\bigl(\CN_x -\sml(\CN_x |_M\smr) \bigr)\;.
\end{equ}
A straightforward calculation, using the symmetries of the problem, shows that property 1 is indeed satisfied.
Property 3 follows immediately from the construction, so it remains to check that property 2 holds, \ie that
\begin{equ}
\CN_3(L_3) \ge 1- b\;,
\end{equ}
for $|b| < 1$, and $\CN_3(L_3)>0$ otherwise. It follows from the definition of the total variation distance $\|\cdot\|_\TV$ that 
\begin{equ}
\CN_3(L_3) = 1-{1 \over 2}\|\sml(\CN_x |_M\smr) - \tau_a^* \sml(\CN_x |_M\smr) \|_\TV\;,
\end{equ}
where $\tau_a(x) = x-a$. 
Since $M\ge 4b \ge 4a$, is clear from the picture and from the fact that the density of the normal distribution is
everywhere positive, that $\CN_3(L_3)>0$ for every $a\in\R$. It therefore suffices to consider the case
$|b| \le 1$. Since
$\int_M^\infty e^{-x^2/2}\,dx < b/8$, one 
has $\|\CN_x|_M -\CN_x \|_\TV \le b/4$, which implies
\begin{equ}
\CN_3(L_3) \ge 1- {b\over 4} - {1 \over 2}\|\CN_x - \tau_a^* \CN_x\|_\TV\;.
\end{equ}
A straightforward computation shows that, for $|a|\le 1$,
\begin{equ}
\|\CN_x - \tau_a^* \CN_x\|_\TV \le \sqrt{e^{a^2}-1} \le \sqrt{2} a \;,
\end{equ}
and the claim follows.
\end{proof}

We will use the following corollary:

\begin{corollary}\label{cor:Girs}
Let $\Wien$ be the Wiener measure on $\Noise_+$, let $g\in \L^2(\R_+)$ with $\|g\| \le b$, let
$M = \max\{4b,2\log(8/b)\}$, and define
the map $\Psi_g\from\Noise_+\to\Noise_+$ by
\begin{equ}
\sml(\Psi_g w\smr)(t) = w(t) + \int_0^t g(s)\,ds\;.
\end{equ}
Then, there exists a measure $\Wien_{g,b}^2$ on $\Noise_+^2$ such that 
the following properties hold:
\begin{claim}[33]
\item[1.] Both marginals of $\Wien_{g,b}^2$ are equal to the Wiener measure $\Wien$.
\item[2.] If $b \le 1$, one has the bound
\begin{equ}[e:goodset]
\Wien_{g,b}^2 \bigl(\bigl\{(\tilde w_x, \tilde w_y)\,|\, \tilde w_y = \Psi_g(\tilde w_x)\bigr\}\bigr) \ge 1-b\;.
\end{equ}
Furthermore, at fixed $b >0$, the above quantity is always positive and a decreasing function of $\|g\|$.
\item[3.] The set
\begin{equ}
\Bigl\{(\tilde w_x, \tilde w_y)\,\Big|\, \exists \kappa \,:\, \tilde w_y(t) = \tilde w_x(t) + \kappa \int_0^t g(s)\,ds\,,\,
|\kappa| \| g\| \le M\Bigr\}
\end{equ}
has full $\Wien_{g,b}^2$-measure.
\end{claim}
\end{corollary}

\begin{proof}
This is an immediate consequence of the $\L^2$ expansion of white noise, using $g$ as one of the
basis functions and applying \lem{lem:gauss} on that component.
\end{proof}

Given this result (and using the same notations as above), we turn to the construction of the coupling function $\Coupl$ 
for step \step{2}. Given an initial condition $Z = (x_0,y_0,w_x,w_y)$, remember that $g_w$ is defined by
\eref{e:defg_w}. We furthermore define the function $\tilde g_w\from \R_+ \to \R^n$ by
\begin{equ}[e:exprgwtilde2]
\tilde g_w(t) = C \int_{-\infty}^0 {t^{{1\over 2}-H}(-s)^{H-{1\over 2}} \over t-s} g_w(s)\,ds\;,
\end{equ}
with $C$ the constant appearing in \eref{e:formg2}. By \eref{e:formg2}, $\tilde g_w$ is the only
function that ensures that \eref{e:defcoupl} holds  if $\tilde w_x$ and $\tilde w_y$ are related
by \eref{e:defgwtilde}. (Notice that, although \eref{e:exprgwtilde1} seems to differ substantially
from \eref{e:exprgwtilde2}, they do actually define the same function.) Given $Z$ as above and 
$a\in\CA$, denote by $g_{a,Z}$ the restriction of $\tilde g_w$ to the interval $[0,2^N]$ (prolonged
by $0$ outside). It follows from \lem{lem:step1} that there exists a constant $K$ such that if
the coupled process was in an admissible state at the beginning of step \step{1}, then
the a-priori estimate
\begin{equ}[e:aprioristep2]
\|g_{a,Z}\|^2 \equiv \int_0^{2^N}\|g_{a,Z}(s)\|^2\,ds \le C 2^{-2\alpha N} \equiv b_N^2
\end{equ}
holds for some constant $C$.
We thus define $b = \max\{b_N,\|g_{a,Z}\|\}$ and denote by $\Wien_{Z,a}^2$ the restriction of
$\Wien_{g_{a,Z},b}^2$ to the ``good'' set \eref{e:goodset} and by $\tilde\Wien_{Z,a}^2$ its restriction
to the complementary set.

We choose furthermore an arbitrary exponent $\beta$ satisfying the condition
\begin{equ}[e:condbeta]
\beta >  {1\over 1-2\alpha}\;.
\end{equ}
With these notations at hand, we define the coupling function for step \step{2}:
\begin{equ}
T(Z,a) = 2^N\;, \quad \Wien_2(Z,a) = \Wien_{Z,a}^2 \times \delta_{a'} + \tilde \Wien_{Z,a}^2 \times \delta_{a''}\;,
\end{equ}
where
\begin{equ}[e:defstep2]
a' = (2,N+1,\tilde N, 0)\;,\quad a'' = (3,0,\tilde N+1, \tilde t_* 2^{\beta N}\tilde N^{4/(1-2\alpha)})\;,
\end{equ}
for some constant $\tilde t_*$ to be determined in the remainder of this section.
The waiting time in \eref{e:defstep2} has been chosen in such a way that the following holds.

\begin{lemma}\label{lem:step2}
Let $(Z_0,a_0) \in \State^2\times\Noise^2\times\CA$ with $Z_0$ admissible and 
denote by $\CT$ the measure on $\State^2\times\Noise^2$
obtained by evolving it according to the successful realisation of step \step{1}, followed by
$N$ successful realisations of step \step{2}, one failed realisation of step \step{2}, and one waiting
period \step{3}. There exists a constant $\tilde t_*$
such that $\CT$-almost every $Z = (x,y,w_x,w_y)$ satisfies
\begin{equ}
\|x - y\| \le 1+ {1+\C4\over \C2}\;,\quad \CK_\alpha(Z) \le \CK_\alpha(Z_0) + {1 \over 2\tilde N^2}\;,
\end{equ}
where $\tilde N$ denotes the value of the corresponding component of $a_0$.
\end{lemma}

\begin{proof}
We first show the bound on the cost function. Given $Z$ distributed according to $\CT$ as in
the statement, we define $g_w$ by \eref{e:defgw2} as usual. The bounds we get on the function
$g_w$ are schematically depicted in the following figure, where the time interval $[\tilde t_2,t_3]$ 
corresponds to the failed realisation of step \step{2}.
\begin{equ}[e:figBounds]
\mhpastefig{Bounds}
\end{equ}
Notice that, except for the contribution coming from times smaller than $t_1$, we are exactly in the 
situation of \eref{e:defg}.
Since the cost of a function is decreasing under time shifts, the contribution to $\CK_\alpha(Z)$
coming from $(-\infty,t_1]$ is bounded by $\CK_\alpha(Z_0)$. Denote by $g$ the function defined
by
\begin{equ}
g(t) = \cases{g_w(t+t_1) & for $t\in [0, t_3-t_1]$,\cr 0& otherwise.}
\end{equ}
Using the definition of the cost function together with \prop{prop:pastfuture} and the Cauchy-Schwarz 
inequality, we obtain for some
constants $C_1$ and $C_2$ the bound
\begin{equ}
\CK_\alpha(Z) \le \CK_\alpha(Z_0) + C_1 \sqrt{|t_3|^{2H-2}-|t_1|^{2H-2}}\|g\| + C_2 \Bigl|{t_1 \over t_3 - t_1}\Bigr|^{\alpha-{1\over 2}}\|g\|_\alpha\;,
\end{equ}
where $\|\cdot\|$ denotes the $\L^2$-norm. 
Since step \step{1} has length $1$ and the $N$th occurrence of step \step{2} has length $2^{N-1}$,
we have
\begin{equ}
|t_3 - t_1| = 2^{N+1}\;,\quad |t_3| = \tilde t_* 2^{\beta N}\tilde N^{4/(1-2\alpha)}\;.
\end{equ}
In particular, one has $|t_3| > |t_3-t_1|$ if $\tilde t_*$ is larger than $1$.
Since 
\begin{equ}
\sqrt{|t_3|^{2H-2}-|t_1|^{2H-2}} \le |t_3|^{H-{3\over 2}} |t_3 - t_1|^{1\over 2} \le  \Bigl|{t_3 \over t_3 - t_1}\Bigr|^{-{1\over 2}}\;,
\end{equ}
this yields (for a different constant $C_1$) the bound
\begin{equ}
\CK_\alpha(Z) \le \CK_\alpha(Z_0) + C_1  \Bigl|{t_3 \over t_3 - t_1}\Bigr|^{\alpha-{1\over 2}} \|g\|_\alpha
\le\CK_\alpha(Z_0) + C_1 {\tilde t_*^{\alpha -{1\over 2}}2^{-\gamma N}\over \tilde N^2}  \|g\|_\alpha\;,
\end{equ}
where we defined $\gamma = (\beta-1)({1\over 2}-\alpha)$. Notice that \eref{e:condbeta} guarantees that
$\gamma > \alpha$. 

We now bound the $\|\cdot\|_\alpha$-norm of $g$. We know from \lem{lem:step1} that
the contribution coming from the time interval $[t_1,\tilde t_2]$ is bounded by some constant $K$.
Furthermore, by \eref{e:aprioristep2}, we have  for the contribution coming from the interval
$[\tilde t_2, t_3]$ a bound of the type 
\begin{equ}
\int_{\tilde t_2}^{t_3}\|\tilde g(s)\|^2\,ds \le C (N+1)^2\;,
\end{equ}
for some positive constant $C$. This yields for $g$ the bound
\begin{equ}
\|g\|_\alpha \le C (N+1)2^{\alpha N}\;,
\end{equ}
for some other constant $C$. Since $\gamma > \alpha$, there exists a constant $C$ such that
\begin{equ}
\CK_\alpha(Z) \le \CK_\alpha(Z_0) + C {\tilde t_*^{\alpha -{1\over 2}}\over \tilde N^2}\;.
\end{equ}
By choosing $\tilde t_*$ sufficiently large, this proves the claim concerning the increase of the total cost. 

It remains to show that, at the end of step \step{3}, the two realisations of \eref{e:mainequ} 
didn't drift to far apart. 
Define $g_b$ by \eref{e:defgb2} as usual and notice that, by construction, $x_t = y_t$ for $t = \tilde t_2$. 
Writing as before $\rho_t = y_t - x_t$, one has for $t > \tilde t_2$ the estimate
\begin{equ}[e:defrhostep2]
\|\rho_t\| \le {\C4 \over \C2} + \int_{\tilde t_2}^t e^{-\C2(t-s)} \|g_b(s)\|\,ds\;.
\end{equ}
We first estimate the contribution coming from the time interval $[\tilde t_2, t_3]$. Denote by
$\tilde g\from [\tilde t_2, t_3] \to \R^n$ the value $g_w$ would have taken, had the last occurence of step
\step{2} succeeded and not failed
(this corresponds to the dashed curve in \eref{e:figBounds}). 
Defining $\hat g = g_w - \tilde g$,  we have by \eref{e:relG5} that, on the interval $t\in [\tilde t_2,t_3]$,
\begin{equ}[e:defgbstep2]
g_b(t) = C_1 {\hat g(\tilde t_2) \over (t-\tilde t_2)^{{1\over 2}-H}} + C_2 \int_{\tilde t_2}^t {{d\hat g\over ds}(s) \over (t-s)^{{1\over 2}-H}}\,ds\;.
\end{equ}
By Corollary~\ref{cor:Girs} and the construction of the coupling function, 
$\hat g$ is proportional to $g_w$ and, by \eref{e:aprioristep2}, 
we also have for $\hat g$ a bound of the type $\|\hat g\| \le C (N+1)$ (the norm is the $\L^2$-norm over the
time interval $[\tilde t_2, t_3]$). Furthermore, \eref{e:boundderg} yields  $\|{d\hat g\over ds}\| \le C (N+1) 2^{-N}$.
Recall that every differentiable function defined on an interval of length $L$ satisfies
\begin{equ}
|f(t)| \le {\|f\| \over \sqrt L} + \Bigl\|{df \over dt}\Bigr\|\sqrt L \;.
\end{equ}
(The norms are $\L^2$-norms.) Using this to bound
the first term in \eref{e:defgbstep2} and the Cauchy-Schwarz inequality for the second term,
 we get a constant $C$ such that $g_b$ is bounded by
\begin{equ}
\|g_b(t)\| \le C (N+1) \bigl(1 + 2^{-{N\over 2}} (t-\tilde t_2)^{H-{1\over 2}}\bigr)\;.
\end{equ}
From this and \eref{e:defrhostep2}, we get an other constant $C$ such that $\|\rho_t\| \le C (N+1)$ at the
time $t = t_3$. We finally turn to the interval $[t_3,0]$. It follows from \eref{e:relG3} that, for some constant
$C$, we have
\begin{equ}
\|g_b(t)\| \le {1\over 2} + C |t-t_3|^{H-1} \|g\|\;,
\end{equ}
where the term ${1\over 2}$ is the contribution from the times smaller than $t_1$. Since we know by
\eref{e:aprioristep2} and Corollary~\ref{cor:Girs} that the $\L^2$-norm of $g$ is bounded by $C(N+1)$ for some constant
$C$,  we obtain the required estimate by choosing $\tilde t_*$ sufficiently large.
\end{proof}

\begin{remark}\label{rem:step2}
To summarise this subsection, we have shown the following, assuming that the coupled system
was in an admissible state before performing step \step{1} and that step \step{1} succeeded:
\begin{claim}[AA]
\item[1.]  There exists constants $\delta' \in (0,1)$ and $K > 0$ such that the $N$th consecutive
occurrence of step \step{2} succeeds with probability larger than $\max\{\delta', 1-K2^{-\alpha N}\}$.
This occurrence has length $2^{N-1}$.
\item[2.] If the $N$th occurrence of step \step{2} fails and the waiting time for step \step{3} is chosen longer
than $\tilde t_* 2^{\beta N}\tilde N^{4/(1-2\alpha)}$, then the state of the coupled system is again admissible after the end of step \step{3}, provided that 
the  cost $\CK_\alpha(Z)$ at the beginning of step \step{1} was smaller than $1 - {1\over 2 \tilde N^2}$.
\item[3.] The increase in the cost given between the beginning of step \step{1} and the end
of step \step{3} is smaller than ${1 \over 2\tilde N^2}$.
\end{claim}
\end{remark}

Now that the construction of the coupling function $\Coupl$ is completed,
we can finally turn to the proof of the results announced in the introduction.

\section{Proof of the main result}
\label{sec:mainproof}

Let us first reformulate Theorems~\ref{theo:smallH} and \ref{theo:largeH} in a more precise way, 
using the notations developed in this paper.

\begin{theorem}\label{theo:mainTheoFormal}
Let $H \in (0,1) \setminus \{{1\over 2}\}$, let $f$ and $\sigma$ satisfy \textbf{A1}--\textbf{A3} if $H<{1\over 2}$ and
\textbf{A1}, \textbf{A2'}, \textbf{A3} if $H>{1\over 2}$, and let $\gamma < \max_{\alpha < H} \alpha(1-2\alpha)$.
Then, the SDS defined in Proposition~\ref{prop:defSDS} has a unique invariant measure $\mu_*$.
Furthermore, there exist positive constants $C$ and $\delta$ such that, for every generalised initial 
condition $\mu$, one has
\begin{equ}[e:convergence]
\|\Evol \CQ_t \mu - \Evol \mu_*\|_\TV \le 2\mu\bigl(\{\|x_0\| > e^{\delta t}\}\bigr) + C t^{-\gamma}\;.
\end{equ}
\end{theorem}

\begin{proof}
The existence of $\mu_*$ follows from \prop{prop:Lyap} and \lem{lem:Lyap}. Furthermore, the assumptions
of \prop{prop:faithful} hold by the invertibility of $\sigma$, so the uniqueness of $\mu_*$ will follow
from \eref{e:convergence}.

Denote by $\phi$ the SDS constructed in Proposition~\ref{prop:defSDS}, and consider the self-coup\-ling 
$\Evol(\mu,\mu_*)$ for $\phi$ constructed in Section~\ref{sec:defCoupl}. We denote by $(x_t,y_t)$
the canonical process associated to $\Evol(\mu,\mu_*)$ and we define a random time $\tilde \tau_\infty$ by
\begin{equ}
\tilde \tau_\infty = \inf \bigl\{t > 0 \,|\, x_s = y_s\, \forall \, s \ge t\bigr\}\;.
\end{equ}
It then follows immediately from \eref{e:TVprop} that
\begin{equ}
\|\Evol \CQ_t \mu - \Evol \mu_*\|_\TV \le 2\prob \sml(\tilde \tau_\infty > t\smr)\;.
\end{equ}
Remember that $\Evol(\mu,\mu_*)$ was constructed as the marginal of the law of a Markov
process with continuous time, living on an augmented phase space $\BigSpace$.
Since we are only interested
in bounds on the random time $\tilde \tau_\infty$ and since we know that $x_s = y_s$ as long
as the coupled system is  in the state \step{2}, it suffices to consider the Markov chain $(Z_n,\tau_n)$
constructed in \eref{e:transZn}. It is clear that $\tilde \tau_\infty$ is then dominated by the random time
$\tau_\infty$ defined as
\begin{equ}
\tau_\infty = \inf \bigl\{\tau_n\,|\, S_m = 2\,\forall\, m\ge n\}\;,
\end{equ}
where $S_n$ is the component of $Z_n$ indicating the type of the corresponding step.
Our interest therefore only goes to the dynamic of $\tau_n$ and $S_n$. We define the
sequence of times $t(n)$ by
\begin{equ}[e:deftn]
t(0) = 1\;,\quad t(n+1) = \inf\{m > t(n)\,|\, S_m = 1\}\;,
\end{equ}
and the sequence of durations $\Delta\tau_n$ by
\begin{equ}
\Delta\tau_n = \tau_{t(n+1)} - \tau_{t(n)}\;,
\end{equ}
with the convention $\Delta\tau_n = +\infty$ if $t(n)$ is infinite (\ie if the set in \eref{e:deftn} is empty).
Notice that we set $t(0) = 1$ and not $0$ because we will treat step \step{0} of the coupled process
separately.
The duration $\Delta\tau_n$ therefore measures the time needed by the coupled system starting
in step \step{1} to come back again to step \step{1}. We define the sequence $\xi_n$ by
\begin{equ}
\xi_0 = 0\;,\quad \xi_{n+1} = \cases{-\infty & if $\Delta\tau_n = +\infty$,\cr \xi_n + \Delta\tau_n& otherwise.}
\end{equ}
By construction, one has
\begin{equ}[e:proptauC]
\tau_\infty = \tau_1 + \sup_{n\ge 0}\xi_n\;,
\end{equ}
so we study the tail distribution of the $\Delta\tau_n$.

For the moment, we leave the value $\alpha$ appearing throughout the paper free, we will tune it at the
end of the proof. Notice also that, by Remarks~\ref{rem:step1} and \ref{rem:step2}, the cost increases by
less than ${1 \over 2\tilde N^2}$ every time the counter $\tilde N$ is increased by $1$. Since the initial
condition has no cost (by the choice \eref{e:definit} of its distribution), this implies that, with probability
$1$, the system is in an admissible state every time step \step{1} is performed.

Let us first consider the probability of $\Delta\tau_n$ being infinite. By \rem{rem:step1}, the probability for
step \step{1} to succeed is always greater than $\delta$. After step \step{1}, the $N$th occurrence of step \step{2}
has length $2^{N-1}$, and a probability greater than $\max\{\delta', 1-K2^{-\alpha N}\}$ of succeeding.
Therefore, one has
\begin{equ}
\prob\sml(\Delta\tau_n \ge 2^N\smr) \ge \delta \prod_{k=0}^N \max\{\delta', 1-K2^{-\alpha k}\}\;.
\end{equ}
This product always converges, so there exists a constant $p_* > 0$ such that
\begin{equ}
\prob\sml(\Delta\tau_n = \infty\smr) \ge p_*\;,
\end{equ}
for every $n > 0$. Since our estimates are uniform over all admissible initial conditions and
the coupling is chosen in such a way that the system is always in an admissible state at the beginning
of step \step{1}, we actually just proved that the conditional probability of $\prob\sml(\Delta\tau_n = \infty\smr)$
on any event involving $S_m$ and $\Delta\tau_m$ for $m < n$ is bounded from below by $p_*$.

For $\Delta\tau_n$ to be finite, there has to be a failure of step \step{2} at some point (see \eref{e:algo}).
Recall that if step \step{2} succeeds exactly $N$ times, the corresponding value for $\Delta\tau_n$ will
be equal to $2^N + \tilde t_* 2^{\beta N} (1+n)^{4/(1-2\alpha)}$ for $N>0$ and to $t_* (1+n)^{4/(1-2\alpha)}$ for
$N=0$. This follows from \eref{e:defa2} and \eref{e:defstep2}, noticing that $\tilde N$ in those formulae
counts the number of times step \step{1} occurred and is therefore equal to $n$. We also know that
the probability of the $N$th occurrence of step \step{2} to fail is bounded from above by $K2^{-\alpha N}$.
Therefore, a very crude estimate yields a constant $C$ such that
\begin{equ}
\prob\bigl((1+n)^{-4/(1-2\alpha)} \Delta\tau_n \ge C2^{\beta N} \,\,\text{and}\,\,\Delta\tau_n \neq \infty\bigr) \le K\sum_{k>N}2^{-\alpha k}\;.
\end{equ}
This immediately yields for some other constant $C$
\begin{equ}[e:boundDeltaT]
\prob\bigl((1+n)^{-4/(1-2\alpha)} \Delta\tau_n \ge T \,\,\text{and}\,\,\Delta\tau_n \neq \infty\bigr) \le C T^{-\alpha/\beta}\;.
\end{equ}
As a consequence, the process $\xi_n$ is stochastically dominated by the Markov chain
$\zeta_n$ defined by
\begin{equ}
\zeta_0 = 0\;,\quad \zeta_{n+1} = \cases{-\infty & with probability $p_*$,\cr\zeta_n +  (n+1)^{4/(1-2\alpha)}p_n & 
with probability $1-p_*$,}
\end{equ}
where the $p_n$ are positive i.i.d.\ random variables with tail distribution $C T^{-\alpha/\beta}$, \ie
\begin{equ}
\prob(p_n \ge T) = \cases{C T^{-\alpha/\beta} & if $C T^{-\alpha/\beta} < 1$, \cr 1 & otherwise.}
\end{equ}
With these notations and using the representation \eref{e:proptauC}, $\tau_\infty$ is bounded by
\begin{equ}[e:boundtauC]
\prob \sml(\tau_\infty > t\smr) \le \prob(\tau_1 > t/2) + \prob \Bigl(\sum_{n=0}^{n_*} (n+1)^{4/(1-2\alpha)} p_n > t/2\Bigr)\;,
\end{equ}
where $n_*$ is a random variable independent of the $p_n$ and such that
\begin{equ}[e:boundnstar]
\prob(n_* = k) = p_*(1-p_*)^k\;.
\end{equ}
In order to bound the second term in \eref{e:boundtauC}, it thus suffices to estimate terms of the form
$\sum_{n=0}^{k} (n+1)^{4/(1-2\alpha)} p_n$ for fixed values of $k$.
Using the Cauchy-Schwartz inequality, one obtains the existence of positive constants $C$ 
and $N$ such that
\begin{equ}
 \prob \Bigl(\sum_{n=0}^{k} (n+1)^{4/(1-2\alpha)} p_n > t/2\Bigr) \le C (k+1)^{N} t^{-\alpha/\beta}\;.
\end{equ}
Combining this with \eref{e:boundnstar} and \eref{e:boundtauC} yields, for some other constant $C$,
\begin{equ}
\prob \sml(\tau_\infty > t\smr) \le \prob(\tau_1 > t/2) + C t^{-\alpha/\beta}\;.
\end{equ}
By the definition of step \step{0} \eref{e:defCoupl0}, we get for $\tau_1$:
\begin{equ}
\prob(\tau_1 > t/2) \le \mu\bigl(\{\|x_0\| > e^{\C2 t/2}/2\}\bigr) + \mu_*\bigl(\{\|y_0\| > e^{\C2 t/2}/2\}\bigr)\;.
\end{equ}
Since, by \prop{prop:Lyap}, the invariant measure $\mu_*$ has bounded moments, the second
term decays exponentially fast. Since $\alpha < \min\{{1\over 2}, H\}$ and $\beta > (1-2\alpha)^{-1}$
are arbitrary, one can realise $\gamma = \alpha/\beta$ for $\gamma$ as in the statement.

This concludes the proof of \theo{theo:mainTheoFormal}.
\end{proof}




We conclude this paper by discussing several possible extensions of our result. The first two extensions
are straightforward and can be obtained by simply rereading the paper carefully and (in the
second case) combining its results with the ones obtained in the references. 
The two other extensions are less obvious and merit further investigation.

\subsection{Noise with multiple scalings} 

One can consider the case where the equation is driven by several independent fBm's with different values of
the Hurst parameter:
\begin{equ}
dx_t = f(x_t)\,dt + \sum_{i=1}^m \sigma_i \, dB^i_{H_i}(t)\;.
\end{equ}
It can be seen that in this case, the invertibility of $\sigma$ should be replaced by the condition that
the linear operator
\begin{equ}
\sigma = \sigma_1 \oplus \sigma_2 \oplus \ldots \oplus \sigma_m \from \R^{mn}\to \R^n\;,
\end{equ}
has rank $n$. The condition on the convergence exponent $\gamma$ then becomes
\begin{equ}
\gamma < \min\sml\{{\gamma_1,\ldots,\gamma_m}\smr\}\;,
\end{equ}
where $\gamma_i = \max_{\alpha < H_i} \alpha(1-2\alpha)$.

\subsection{Infinite-dimensional case}

In the case where the phase space for \eref{e:mainequ} is infinite-dimensional, the question of 
global existence of solutions is technically more involved and was tackled in \cite{MaslFBM}. 
Another technical difficulty arises from the fact that one might want to take for $\sigma$ an
operator which is not boundedly invertible, so \textbf{A3} would fail on a formal level.
One expects to be able to overcome this difficulty at least in the case where the equation is semilinear
and parabolic, \ie of the type
\begin{equ}
dx = Ax\,dt + F(x)\,dt + Q\,dB_H(t)\;,
\end{equ}
with the domain of $F$ ``larger'' (in a sense to be quantified) than the domain of $A$
and $B_H$ a cylindrical fBm on some Hilbert space $\CH$ on which the solution
is defined, provided the eigenvalues of $A$ and of $Q$ satisfy some compatibility
condition as in \cite{ZDP1,Cerr,EH3}.

On the other hand, it is possible in many cases to split the phase space into
a finite number of ``unstable modes'' and an infinite number of ``stable modes'' that
are slaved to the unstable ones. In this situation, it is sufficient to construct step \step{1}
in such a way that the unstable modes meet, since the stable ones will then automatically
converge towards each other. A slight drawback of this method is that the convergence
towards the stationary state no longer takes place in the total variation distance.
We refer to \cite{MatNS,ArmenKuk,HExp02} for implementations
of this idea in the Markovian case.

\subsection{Multiplicative noise}

In this case, the problem of existence of global solutions can already be hard. In the case $H>1/2$, the
fBm is sufficiently regular, so one obtains pathwise existence of solutions by rewriting \eref{e:mainequ} in
integral form and interpreting the stochastic integral pathwise as a Riemann-Stieltjes integral.
In the case $H\in ({1\over 4}, {1\over 2})$,
it has been shown recently \cite{MR96b:60150,MR2000c:60089,MR2003c:60066} that pathwise 
solutions can also be obtained by realising the fBm as
a geometric rough path. More refined probabilistic estimates are required in the analysis
of step \step{1} of our coupling construction. The equivalent of equation \eref{e:equrho} 
then indeed contains a multiplicative noise term, so the deterministic estimate \eref{e:estrho} fails.

\subsection{Arbitrary Gaussian noise}

Formally, white noise is a centred Gaussian process $\xi$ with correlation function
\begin{equ}
\expect \xi(s)\xi(t) = C_w(t-s) = \delta(t-s)\;.
\end{equ}
The derivative of the fractional Brownian motion with Hurst parameter $H$ is formally also a centred Gaussian
process, but its correlation function is proportional to
\begin{equ}
C_H(t-s) = |t-s|^{2H-2}\;,
\end{equ}
which should actually be interpreted as the second derivative of $|t-s|^{2H}$ in the sense of distributions.

A natural question is whether the results of the present paper also apply to differential equations driven by
Gaussian noise with an arbitrary correlation function $C(t-s)$. There is no conceptual obstruction to
the use of the method of proof presented in this paper in that situation, but new estimates are required.
It relies on the fact that the driving process is a fractional Brownian motion only to be able to explicitly
perform the computations of Section~\ref{sec:defCoupl}. 

\bibliographystyle{Martin}
\markboth{\sc \refname}{\sc \refname}
\bibliography{refs}
\end{document}